\documentclass[12 pt]{article}
\counterwithin*{equation}{section}
\usepackage[all]{xy}
\usepackage{amsmath,amssymb,amsfonts,latexsym,float,graphics,epsfig}

\usepackage{setspace}
\usepackage{appendix}
\usepackage{xcolor}
\usepackage[top=2cm, bottom=2cm, left=2cm, right=2cm]{geometry}
\usepackage[colorlinks]{hyperref}
\usepackage{tikz}
\usetikzlibrary{cd}

\newtheorem{theorem}{Theorem}[section]
\newtheorem{proposition}[theorem]{Proposition}

\begin{document}

\title{On Locally Conformally Cosymplectic Hamiltonian Dynamics and Hamilton-Jacobi Theory}
\author{}
\maketitle
\begin{center}

Begüm Ateşli \footnote{E-mail: 
\href{mailto:b.atesli@gtu.edu.tr}{b.atesli@gtu.edu.tr}, corresponding author}\\
Department of Mathematics, \\ Gebze Technical University, 41400 Gebze,
Kocaeli, Turkey.

\bigskip

O\u{g}ul Esen\footnote{E-mail: 
\href{mailto:oesen@gtu.edu.tr}{oesen@gtu.edu.tr}}\\
Department of Mathematics, \\ Gebze Technical University, 41400 Gebze,
Kocaeli, Turkey.

\bigskip

Manuel de Le\'on\footnote{E-mail: \href{mailto:mdeleon@icmat.es}{mdeleon@icmat.es}}
\\ Instituto de Ciencias Matem\'aticas, Campus Cantoblanco \\
 Consejo Superior de Investigaciones Cient\'ificas
 \\
C/ Nicol\'as Cabrera, 13--15, 28049, Madrid, Spain
\\
and
\\
Real Academia Espa{\~n}ola de las Ciencias.
\\
C/ Valverde, 22, 28004 Madrid, Spain.

\bigskip 

Cristina Sard\'on\footnote{E-mail: \href{mailto:mariacristina.sardon@upm.es}{mariacristina.sardon@upm.es}}
\\ Department of Applied Mathematics 
\\ Universidad Polit\'ecnica de Madrid 
\\ C/ Jos\'e Guti\'errez Abascal, 2, 28006, Madrid. Spain.

\end{center}

\begin{abstract}
Cosymplectic geometry has been proven to be a very useful geometric background to describe time-dependent Hamiltonian dynamics. In this work, we address the globalization problem of locally cosymplectic Hamiltonian dynamics that failed to be globally defined. We investigate both the geometry of locally conformally cosymplectic (abbreviated as LCC) manifolds and the Hamiltonian dynamics constructed on such LCC manifolds. Further, we provide a geometric Hamilton-Jacobi theory on this geometric framework. 
\end{abstract}

\tableofcontents

\setlength{\parskip}{4mm}

\onehalfspacing
\section{Introduction}

As it is very well known, symplectic manifolds provide a convenient geometrical frameworks for classical autonomous Hamiltonian dynamics. We cite an incomplete list of works \cite{AbrahamMarsden,Arnold-book,leon89,holm2009geometric} setting the fundamentals of geometric Hamiltonian dynamics.

Symplectic manifolds are even dimensional, therefore, if one wants to consider time on the framework, i.e., we are in the case of $t$-dependent Hamiltonian dynamics, one needs to add the time parameter as a new variable. This implies that we are working on an odd dimensional manifold, which is beyond symplectic geometry.  The odd dimensional counterpart of symplectic geometry is cosymplectic geometry \cite{Libermann55,Liber87,Takizawa63}. As pointed out in the literature, cosymplectic manifolds provide a proper geometrical framework for time-dependent (non-autonomous) Hamiltonian dynamical systems. We also cite a recent review paper \cite{cape} covering this issue and more on cosymplectic manifolds. 

In this work, our aim is to investigate the gluing problem of local cosymplectic Hamiltonian flows and study the associated Hamilton-Jacobi formalism. This provides a generalization of cosymplectic dynamics \cite{LeonSolanoMarrero96,GuzmMarr10,
LaciMarrPadr12} as well as the geometric Hamilton-Jacobi theory on cosymplectic manifolds \cite{LeonSar2}. To be more precise about our goal and to state the problem more concretely, we provide the following discussions.

\textbf{Almost Cosymplectic Manifolds.}
A $2n+1$ dimensional almost cosymplectic manifold admits a differential one-form $\eta$ and a differential two-form $\Omega$ so that the top-form $\eta\wedge\Omega^n\neq 0$. This geometry permits us to define an almost Poisson bracket (a bracket satisfying the Leibnitz identity but not necessarily the Jacobi identity) on the space of real valued functions. Actually, almost cosymplectic manifolds are particular instances of Leibniz manifolds \cite{OrtePlan04}. An almost cosymplectic manifold is called cosymplectic if $\eta$ and $\Omega$ are closed forms. In this case, the induced bracket satisfies the Jacobi identity and turns out to be a Poisson bracket \cite{Laurent13,vaisman2012lectures,Weinstein98}. An interesting subclass of almost cosymplectic     
manifolds are locally conformally cosymplectic manifolds.



\textbf{Locally Conformally Cosymplectic Manifolds.}
Our approach to LCC manifolds is as follows. We start with an almost cosymplectic manifold equipped with $(\eta,\Omega)$. Evidently, in each local chart this manifold remains to be almost cosymplectic. Locally, a LCC manifold admits an atlas whose charts are cosymplectic. These local cosymplectic structures do not necessarily glue up to a global one, this is why a LCC manifold fails to be a cosymplectic manifold, but the globalization problem of LCC can be achieved by means of some local conformal parameters. In short, we can say that LCC manifolds are almost cosymplectic manifolds that locally behave as cosymplectic manifolds up to some conformal factors. A similar discussion is also available for the case of symplectic geometry, in which one defines locally conformally symplectic (abbreviated as LCS) manifolds as almost symplectic manifolds with local symplectic character up to conformal factors, see \cite{Bazzoni2018,Hwa,Vaisman85}. 

In this work, we are interested in the Hamiltonian formulation on LCC manifolds. The situation is the following.
In local charts, referring to the local cosymplectic formalism, we write cosymplectic Hamiltonian dynamics. Then, by means of the conformal parameters, we glue the local dynamics to a global one. In this regard, we may state from the physical point of view that LCC manifolds provide a proper geometric background for gluing all of local time-dependent Hamiltonian dynamics \cite{LeonLCCos}.  The new results included in this work are listed in the next lines. 

\textbf{(1) ''Symplectization" of LCC Manifolds.} 
In \cite{Leon-CosympReduction}, it is shown that there exists a symplectic form on the  trivial line bundle over the cosymplectic bundle. We may call this symplectization of cosymplectic manifolds, following the terminology of Arnold \cite{Arnold-book}. Accordingly, we shall state in Proposition \ref{LCCos-to-LCS} that a line bundle over a LCC manifold is a LCS manifold. 

\textbf{(2) Darboux Theorem.} Darboux theorems are important to characterize manifolds equipped with certain tensorial objects. It is known that cotangent bundles admit LCS structures \cite{HaRy99A,ChanMurp19}. In \cite{OtSt15}, it is shown that locally every LCS manifold has this form. There exists a Darboux theorem on cosymplectic manifolds determining the generic example as the extended cotangent bundle, see \cite{cape}.   In the light of these observations, we shall establish Proposition \ref{Darboux-LCC} determining a Darboux theorem for LCC manifolds.

\textbf{(3) Jacobi Structure of LCC Manifolds.} LCS manifolds are Jacobi  manifolds, we shall explicitly show the Jacobi character of LCC manifolds in Proposition \ref{prop-LCC-Jacobi}. 

\textbf{(4) Algebra of One-Forms: Lie Algebroid Realizations.} It is known that Jacobi manifolds admit Lie algebroid formalisms \cite{LeonMarr97}. We shall review this for LCS manifolds in Section \ref{Sec-LCS-Algebroid} and show that Lichnerowicz-deRham (abbreviated as LdR) exact one-forms constitute a subalgebra. In Section \ref{Sec-LCC-Algebroid} we shall write an algebra on the space of one-form sections on LCC. Then we shall discuss the Lie algebroid realization of LCC manifolds. 

\textbf{(5) Hamilton-Jacobi Formalism for LCC Hamiltonian Dynamics.} 
Hamilton-Jacobi theory provides a way to solve Hamilton equations. 
In a recent paper \cite{EsenLeonSarZaj1}, we  examined the geometric Hamilton-Jacobi theory for Hamiltonian dynamics on LCS manifolds. Additionally, we took the locally conformal discussions to $k$-symplectic formalism \cite{EsLeSaZa-k-sympl} and jet bundle formalism \cite{EsLeSaZa-Cauchy} in order to define locally conformally Hamiltonian field theories. In that work, we provided the locally conformal Hamilton-De Donder-Weyl  formalism, as well as the geometric Hamilton-Jacobi theories for these extensions. In this work, we shall present two versions (namely in Theorem \ref{gamma-rel-thm-LCCos} and Theorem \ref{gamma-rel-thm-LCC}) of the geometric Hamilton-Jacobi theorem in the realm of LCC Hamiltonian dynamics. These results are a generalization of the HJ theorems obtained for cosymplectic Hamiltonian dynamics in \cite{LeonSar2}. 

\textbf{The content.} In the following section we summarize Hamiltonian dynamics on symplectic and LCS manifolds and their corresponding geometric HJ theorems. In Section \ref{Sec-Cos}, we shall recall cosymplectic manifolds and the Hamilton-Jacobi theorem for cosymplectic Hamiltonian dynamics. Then we shall explain the basics on LCC geometry, and it is in this same section where the symplectization problem of LCC manifolds and Darboux coordinates of LCC are established. 
In Section \ref{Sec-Dyn-LCC-Ev}, dynamics on LCC manifolds is explained, Jacobi and Lie algebroid characters of LCC manifolds will be obtained, and finally the Hamilton-Jacobi theorems for LCC Hamiltonian dynamics will be written.

\textbf{Notation.} From now on we consider $\mathfrak{X}(M)$ to be the space of vector fields on $M$, whereas $\Gamma^k(M)$ is the space of $k$-form sections on $M$. $\mathcal{L}_{X}$ is the Lie derivative with respect to the vector field $X$. An arbitrary almost symplectic manifold is represented by the pair $(N,\omega)$ while an arbitrary almost cosymplectic manifold is denoted by $(M,\eta,\Omega)$. On $T^*Q$ we consider the differential forms $\theta_Q$ and $\omega_Q=-d\theta_Q$ and their pull backs to $T^*Q\times \mathbb{R}$ are $\Theta_Q$ and $\Omega_Q$, respectively.

\textbf{Abbreviations.} 
We refer frequently the following  abbreviations: 
\begin{center}
\begin{tabular}{||c | c ||} 
 \hline
 Abbreviation &  Extension  \\ [0.5ex] 
 \hline \hline
 HJ & Hamilton-Jacobi  \\ 
 \hline
 LdR &  Lichnerowicz-deRham   \\
 \hline
 LCS &  locally conformal symplectic \\
 \hline
 LCC & locally conformal cosymplectic   \\
 \hline
\end{tabular}
\end{center}

\section{Locally Conformal Symplectic Hamiltonian Dynamics}

\subsection{Symplectic Dynamics and HJ Theory}

A manifold $N$ equipped with a non-degenerate two-form $\omega$ is called an almost symplectic manifold \cite{vaisman13}. The non-degeneracy of the two-form manifests that $N$ is necessarily even dimensional. The musical mapping 
\begin{equation}
\omega^{\flat }:\mathfrak{X}\left(N\right)
\longrightarrow \Gamma^{1} (N),\qquad X\mapsto \iota_{X}\omega
_{N},\label{bemol}
\end{equation}
induced from the two-form $\omega$ is an isomorphism with the inverse mapping $\omega^\sharp$. Here, $\iota_{X}$ is the interior derivative (contraction) operator. 
An almost symplectic manifold is called symplectic if $\omega$ is closed as well. We cite an incomplete list of references \cite{AbrahamMarsden,Arnold-book,leon89,holm2009geometric,Liber87}.  

For a Hamiltonian function $H$, on a symplectic manifold $(N,\omega) $, a Hamiltonian vector field $X_{H}$  is defined to be
\begin{equation}
\iota_{X_{H}}\omega =dH.  \label{Hamvf}
\end{equation}
for a real valued function $H$ on $N$ (the function $H$ is the Hamiltonian). Here, $\iota_ {X_H}$ denotes the interior derivative with respect to $X_H$.

\textbf{Cotangent Bundle.} The cotangent bundle $T^{\ast }Q$ of a manifold $Q$ is a canonically symplectic manifold admitting the Liouville one-form $\theta_Q$ and the canonically symplectic two-form $\omega_Q=-d\theta_Q$. The value of $\theta_Q$ over a vector field $X$ on $T^{\ast }Q$  is defined to be 
\begin{equation}
\theta_{Q}( X) =\left\langle \tau
_{T^{\ast }Q} ( X  ) ,T\pi _{Q} ( X
 ) \right\rangle ,  \label{canonicaloneform}
\end{equation}
where $\tau_{T^{\ast }Q}:TT^*Q\mapsto T^*Q$ is the tangent bundle projection whereas $T\pi _{Q}$ is the tangent mapping of the cotangent bundle projection $\pi_Q:T^*Q\mapsto Q$. 
 
 \textbf{HJ for Symplectic Dynamics.} Consider the Hamiltonian vector field $X_H$  defined on the symplectic manifold $T^*Q$. Define a projection of the Hamiltonian vector field by means of a one-form section $\varrho$ on $Q$ as follows
  \begin{equation}\label{X-H-proj}
 X_H^{\varrho}= T\pi_Q \circ X_H \circ \varrho.
 \end{equation}
 See that $X_H^{\varrho}$ is a vector field on $Q$. 
We provide the following diagram to visualize this projection  
 \begin{equation}\label{geomdiagram--}
\xymatrix{ T^{*}Q
\ar[dd]^{\pi} \ar[rrr]^{X_H}&   & &TT^{*}Q\ar[dd]^{T\pi}\\
  &  & &\\
 Q\ar@/^2pc/[uu]^{\varrho}\ar[rrr]^{X_H^{\varrho}}&  & & TQ}
\end{equation}

The geometric Hamilton-Jacobi theory relates the projected vector field $X_H^{\varrho}$ and the Hamiltonian vector field $X_H$. We refer to  \cite{CariHJ}  for further explanations. 

\begin{theorem}\label{th1}
If $\varrho$ is a closed one-form, then the following conditions are equivalent:
\begin{enumerate}
\item The vector fields $X_H$ and $X_H^{\varrho}$ are $\varrho$-related that is $T\varrho \circ X_H^{\varrho} = X_H  \circ \varrho$. 
\item The equation $d(H\circ \varrho)=0$ holds.
\end{enumerate}
\end{theorem}
The $\varrho$-relatedness stated in the first item gives that if $q(t)$
is an integral curve of the projected field $X_{H}^{\varrho}$, then $\varrho\circ q(t)$
is an integral curve of the Hamiltonian vector field $X_{H}$. Thanks to the Poincar\'{e} lemma, the second condition gives that, there locally  exists a function $S$ on $Q$ satisfying $dS=\varrho$ so that 
\begin{equation}\label{HJ1}
H\left(q^i,\frac{\partial S}{\partial q^i}\right)=\epsilon,
\end{equation} 
where $\epsilon$ is a constant. This is the classical form of the time independent Hamilton-Jacobi equation \cite{Arnold-book,Goldstein-book}. We refer to a recent survey on HJ \cite{esen2022reviewing} for some more details on the present discussion.

\subsection{Locally Conformal Symplectic Dynamics and HJ Theory}\label{Sec-LCS}

An almost symplectic manifold $(N,\omega)$ is called locally conformal symplectic (abbreviated as LCS) manifold if the two-form $\omega$ is closed locally up to a conformal factor \cite{Bazzoni2018,Hwa,Vaisman85}. A LCS manifold admits an open covering by a family of open sets $U_\alpha$ with an associated set of functions $\rho_{\alpha}$ such that 
\begin{equation}
d(e^{-\rho_{\alpha}}\omega\big\vert_{\alpha})=0
\end{equation}
where $\omega\vert_{\alpha}$ is the restriction of the almost symplectic structure $\omega$ to $U_\alpha$. This permits us to define a local two-form $
    \omega_{\alpha}=e^{-\rho_{\alpha}}\omega\vert_{\alpha}$
which turns out to be symplectic. In the intersection $U_\alpha \cap U_\beta$ of two open neighbourhoods, the local symplectic two-forms are related in the following way
\begin{equation} \label{transition}
\omega_{\beta}=e^{-(\rho_{\beta}-\rho_{\alpha})}\omega_{\alpha},\qquad 
\kappa_{\beta \alpha}=e^{\rho_{\alpha}}/e^{\rho_{\beta}}=e^{-(\rho_{\beta}-\rho_{\alpha})}.
\end{equation}
In the intersection of three open neighbourhoods, 
the scalars satisfy the cocycle condition $\kappa_{\beta \alpha}\kappa_{\alpha \rho}=\kappa_{\beta \rho}$, so that 
we can glue the local symplectic two-forms $\omega_\alpha$ up to a line bundle $L \mapsto N$ valued two-form on $N$. 

\textbf{LCS Manifold with Lee-form.}
It can easily be seen that the exterior derivative of the local functions coincide in the intersection of two neighbourhoods that is $d\rho_{\alpha}=d\rho_{\beta}$. By gluing up these exact one-forms, we arrive at a closed one-form, called the Lee one-form and denoted by $\theta$, that takes the local form $d\rho_{\alpha}$, see \cite{Hwa}. This observation permits us to get an alternative definition of LCS manifold as follows. An almost symplectic manifold $(N,\omega)$ is called a LCS manifold if there exists  $\theta$ such that
\begin{equation}
d\omega=\theta\wedge \omega.
\end{equation} 
In terms of the LdR differential (see Appendix \ref{Sec-LdR} for the rigorous definition), an almost symplectic manifold $(N,\omega)$ is a LCS manifold if and only if 
\begin{equation}
d_\theta\omega=0
\end{equation} 
for a closed one-form $\theta$. We denote a LCS manifold admitting a Lee-form $\theta$ by the triplet $(N,\omega,\theta)$. Notice that a LCS manifold turns out to be a symplectic manifold if the Lee-form is identically zero. For a generalization of this discussion in the presymplectic formalism, we cite \cite{Wade04}.  

\textbf{Hamiltonian Dynamics.} Consider a LCS manifold $(N,\omega,\theta)$, and a local symplectic space $(U_\alpha,\omega_\alpha)$. For a local Hamiltonian function $H_\alpha$ on $U_\alpha$, we write the geometric Hamilton equation by
\begin{equation}\label{geohamalpha}
    \iota_{X_{\alpha}}\omega_{\alpha}=dH_{\alpha},
\end{equation}
where $X_{\alpha}$ is the local Hamiltonian vector field  function associated to this framework. We glue local Hamiltonian functions $H_\alpha$ to a section of the line bundle $L\mapsto N$. We also glue the following local functions 
\begin{equation} \label{glueHamFunc}
H\vert_\alpha=e^{\rho_\alpha}H_\alpha
\end{equation}
those obeying the intersection law $e^{\rho_{\alpha}}H_{\alpha}=e^{\rho_{\beta}}H_{\beta}$ to a real valued function $H$ on $M$. 
In the light of these realizations, we glue the local Hamilton equation \eqref{geohamalpha} to 
\begin{equation}\label{semiglobal}
   \iota_{X_{H}}\omega=d_\theta H,
\end{equation}
where $X_{H}$ is the vector field obtained by gluing all the vector fields $X_\alpha$. See that, $\theta$ is the Lee-form associated with the LCS manifold and $d_\theta$ stands for the LdR differential defined in (\ref{LdR-Diff}). 


\textbf{LCS Structure on the Cotangent Bundle.}
Let $\psi$  be a closed one-form on the base manifold $Q$ and we pull it back to $T^*Q$ by means of the cotangent bundle projection $\pi_Q$. This determines  a closed semi-basic one-form $\theta=\pi_Q^*(\psi)$ on $T^*Q$. Using the LdR differential, we define
\begin{equation} \label{omega_theta}
\omega_\theta=-d_\theta(\theta_Q)=\omega_Q+\theta\wedge \theta_Q
\end{equation}
on the cotangent bundle $T^*Q$. Since $d\omega_\theta=\theta\wedge \omega_\theta$ holds, the triple $(T^*Q,\omega_\theta;\theta)$ is a LCS manifold  \cite{Banyaga,ChMu17,HaRy99,OtSt15}. It is important to note that all LCS manifolds admit such a local picture.

\textbf{Hamilton-Jacobi Theorem.} 
Consider a Hamiltonian vector field $X_H$ on $T^*Q$ determined through \eqref{semiglobal} and let $\varrho$ be a section of the cotangent bundle fibration $\pi_Q:T^*Q\to Q$. We project $X_H$ to a vector field $X^{\varrho}_H$ on $Q$ by means of $\varrho$ referring to the definition \eqref{X-H-proj}. We state the following Hamilton-Jacobi theorem for Hamiltonian dynamics on a LCS setting, see \cite{EsenLeonSarZaj1}.
\begin{theorem} Consider the LCS manifold $(T^*Q,\omega)$ equipped with the Lee-form $\theta$. If $\varrho$ is a one-form satisfying $d_\psi \varrho= 0$, then the following conditions are equivalent:
\begin{enumerate}
\item The vector fields $X_H$ and $X^{\varrho}_H$ are $\varrho$-related, that is, $
T\varrho \circ X^{\varrho}_H = X_H \circ \varrho. $
\item The equation $d_\psi(H \circ \varrho) = 0$ holds. 
\end{enumerate}
\end{theorem}
We present the following table for a brief summary of LCS manifolds.
     
\begin{center}
\begin{tabular}{ ||c|c||c||  }
 \hline
 \multicolumn{3}{||c||}{\textbf{LCS Manifolds}} \\
 \hline \hline
 Global & Local & Relation\\
 \hline
 $\omega\vert_\alpha$ & $\omega_\alpha$ & $\omega_\alpha = e^{-\rho_\alpha}\omega\vert_\alpha$  \\ 
 $H\vert_\alpha$ & $H_\alpha$ & $H_\alpha = e^{-\rho_\alpha}H\vert_\alpha$\\
$X_H\vert_\alpha$ & $X_\alpha$ & $X_H\vert_\alpha = X_\alpha$\\
 \hline
 $\iota_{X_H} \omega = d_\theta H$ & $\iota_{X_\alpha} \omega_\alpha = dH_\alpha$ & $\theta\vert_\alpha = d\rho_\alpha$\\
 \hline
\end{tabular}
\end{center}

Here, the tensor fields in the global column are those can be glued up to classical (real valued) tensor fields whereas in the local column are those can be glued up to line bundle valued  tensor fields.

\subsection{Algebra of One-forms}\label{Sec-LCS-Algebroid} 

Consider a LCS manifold $(N,\omega)$ and the musical isomorphism $\omega^\flat$ defined in \eqref{bemol} together with its inverse $\omega^\sharp$. 
By definition, we are able to write the identity
\begin{equation}\label{pair-id}
\omega(\omega^\sharp(\mu), X) = \langle \mu, X \rangle.
\end{equation}
for an arbitrary vector field $X$ and an arbitrary one-form $\mu$ on $N$.

\textbf{Jacobi Structure of LCS Manifolds.}
Let us consider now a LCS manifold $(N,\omega)$ with a Lee form $\theta$. Let us now concentrate on the particular case of LCS manifolds. The musical mapping $\omega^\sharp$ permits us to define the Lee vector field $Z_\theta$ as 
\begin{equation} \label{Lee-v-f}
Z_\theta:=\omega^\sharp(\theta), \qquad \iota_{Z_\theta}\omega=\theta.
\end{equation}
For the two one-forms $\mu$ and $\nu$ on $N$, one can define a bivector $\Lambda$ on $N$ as 
\begin{equation} \label{gamma}
\Lambda(\mu,\nu)=\omega(\omega^\sharp(\mu), \omega^\sharp(\nu)).
\end{equation} 
By direct calculation, one can show that the pair $(\Lambda,Z_\theta)$ satisfies the conditions in \eqref{ident-Jac}, see \cite{LeonSar2}, so it determines a Jacobi structure on $N$. Accordingly, the Jacobi bracket of two functions is determined through the following formula
\begin{equation}\label{LCS-brac-Jac}
\{F,H\}_{LCS}=\Lambda(dF,dH)+FZ_\theta(H)-HZ_\theta(F)=\omega(X_F,X_G).
\end{equation} 
Notice that, $X_F$ and $X_H$ are the Hamiltonian vector fields defined as in (\ref{semiglobal}). 

\textbf{The Algebra.} Consider the Jacobi-Lie bracket of vector fields on the manifold $N$. Referring to this Lie algebra structure, we define a Lie bracket on the space of one-form sections $\Gamma^1(N)$ as follows 
\begin{equation}
\{\bullet , \bullet \}_{\Gamma^1(N)}: \Gamma(N) \times \Gamma(N) \longrightarrow \Gamma(N), \qquad (\mu,\nu)\mapsto \omega^{\flat}([\omega^\sharp(\mu),\omega^\sharp(\nu)])
\end{equation}
such that the musical isomorphisms $\omega^{\flat}$ and $\omega^{\sharp}$ turn out to be Lie algebra isomorphisms, that is,
\begin{equation}\label{iso-id}
\omega^\sharp\{\mu,\nu\}_{\Gamma^1(N)}  = [\omega^\sharp(\mu),\omega^\sharp(\nu)], \qquad \omega^\flat[X,Y]  = \{X^\flat,Y^\flat\}_{\Gamma^1(N)}.
\end{equation}
In order to explicitly determine the Lie bracket on the space of one-form sections, we apply the LdR differential to the LCS two-form $\omega$ at $(\omega^\sharp(\mu),\omega^\sharp(\nu),X)$. This gives us that
\begin{equation}\label{LDR-omega}
\begin{split}
d_\theta\omega (\omega^\sharp(\mu),\omega^\sharp(\nu),X) &= d\omega (\omega^\sharp(\mu),\omega^\sharp(\nu),X) - \theta\wedge \omega (\omega^\sharp(\mu),\omega^\sharp(\nu),X)\\
& =\omega^\sharp(\mu)\big(\omega(\omega^\sharp(\nu),X)\big) - \omega^\sharp(\nu)\big(\omega(\omega^\sharp(\mu),X)\big) + X\big(\omega (\omega^\sharp(\mu),\omega^\sharp(\nu))\big) \\
&\qquad -\omega\big([\omega^\sharp(\mu),\omega^\sharp(\nu)],X\big) + \omega\big([\omega^\sharp(\mu),X],\omega^\sharp(\nu)\big) - \omega\big([\omega^\sharp(\nu),X],\omega^\sharp(\mu)\big) \\
&\qquad -\theta(\omega^\sharp(\mu))\omega(\omega^\sharp(\nu),X) + \theta(\omega^\sharp(\nu))\omega(\omega^\sharp(\mu),X) - \theta(X)\omega(\omega^\sharp(\mu),\omega^\sharp(\nu)) \\
&= \omega^\sharp(\mu)\big( \langle\nu,X\rangle \big) - \omega^\sharp(\nu)\big( \langle\mu,X\rangle \big) +\langle d\iota_{\omega^\sharp(\nu)}\iota_{\omega^\sharp(\mu)}\omega, X\rangle \\
&\qquad -\langle\{\mu,\nu\}_{\Gamma^1},X\rangle - \langle \nu,[\omega^\sharp(\mu),X]\rangle + \langle \mu,[\omega^\sharp(\nu),X]\rangle\\
&\qquad -\theta(\omega^\sharp(\mu))\langle\nu,X\rangle + \theta(\omega^\sharp(\nu))\langle\mu,X\rangle - \iota_{\omega^\sharp(\nu)}\iota_{\omega^\sharp(\mu)}\omega\langle\theta,X\rangle.
\end{split}
\end{equation}
Notice that we have used the identities \eqref{pair-id} and \eqref{iso-id} in order to write the last equality above. The commutation property of the Lie derivative and the interior derivative imply
\begin{equation}
\omega^\sharp(\mu)\big( \langle\nu,X\rangle \big) -\langle \nu,[\omega^\sharp(\mu),X]\rangle = \langle\mathcal{L}_{\omega^\sharp(\mu)}\nu,X\rangle, \qquad \omega^\sharp(\nu)\big( \langle\mu,X\rangle \big) -\langle \mu,[\omega^\sharp(\nu),X]\rangle = \langle\mathcal{L}_{\omega^\sharp(\nu)}\mu,X\rangle.
\end{equation}
Here, $\mathcal{L}$ stands for the Lie derivative. 
After we substitute these identities in equation \eqref{LDR-omega}, the fact that $\omega$ is LdR closed and the nondegeneracy of the pairing between the vector fields and one-forms allow us to write the following proposition.
\begin{proposition}
For a LCS manifold $(N,\omega)$, the space $\Gamma^1(N)$ of one-forms is a Lie algebra with the following the Lie bracket 
\begin{equation}\label{Lie-LCS}
 \{\mu,\nu\}_{\Gamma^1(N)} = \mathcal{L}_{\omega^\sharp(\mu)}\nu -  \theta(\omega^\sharp(\mu))\nu - \mathcal{L}_{\omega^\sharp(\nu)}\mu + \theta(\omega^\sharp(\nu))\mu + d_\theta\iota_{\omega^\sharp(\nu)}\iota_{\omega^\sharp(\mu)}\omega.
\end{equation}
\end{proposition}
Notice that, for the case of $\theta=0$, that is, the symplectic subcase, the algebra {\it Lie-LCS} henceforth, reduces to the bracket \cite{Fuch82,Kosz85,vaisman2012lectures}
\begin{equation}\label{Lie-LCS-Sym}
 \{\mu,\nu\}_{\Gamma^1(N)} = \mathcal{L}_{\omega^\sharp(\mu)}\nu  - \mathcal{L}_{\omega^\sharp(\nu)}\mu  + d\iota_{\omega^\sharp(\nu)}\iota_{\omega^\sharp(\mu)}\omega.
\end{equation}

\textbf{Lie Algebroid Realization.} As it is well known, the cotangent bundle of a Poisson manifold admits a Lie algebroid realization \cite{CostDazoWein87}. The bracket \eqref{Lie-LCS} generalizes this to a more general setting. We draw the following bundle mapping 
 \begin{equation}\label{geomdiagram-}
\xymatrix{ (T^{*}N,\{\bullet,\bullet\}_{\Gamma^1})
\ar[ddr]_{\pi_N} \ar[rr]^{\quad \omega^\sharp}&   & TN \ar[ddl]^{\tau_N}\\
  &    &\\
 & N}
\end{equation}
Referring to Lie algebroid definition given in Appendix \ref{Defn-Lie-Algebroid}, we state that for a LCS manifold $(N,\omega)$, the quintuple $(T^*N,\pi_N,N,\omega^\sharp,\{\bullet,\bullet\}_{\Gamma^1(N)})$ is a Lie algebroid. Here, the anchor map is the musical isomorphism $\omega^\sharp$ according to the identifications in \eqref{iso-id}, whereas the Lie algebroid bracket is the one in \eqref{Lie-LCS}. It is needless to say that $\theta$ identically vanishes if LCS reduces to a symplectic manifold and that the Lie algebroid structure reduces to the one obtained in the classical Poisson case, see \cite{Weinstein98}. 

Let us show that the identity in \eqref{oid-id} is indeed satisfied. Let $F\in\mathcal{F}(N)$ and $\mu,\nu \in \Gamma^1(N)$. We have
\begin{equation}
    \begin{split}
        \{\mu,F\nu\}_{\Gamma^1(N)} &= \mathcal{L}_{\omega^\sharp(\mu)}(F\nu) -  \theta(\omega^\sharp(\mu))F\nu - \mathcal{L}_{\omega^\sharp(F\nu)}\mu + \theta(\omega^\sharp(F\nu))\mu + d_\theta\iota_{\omega^\sharp(F\nu)}\iota_{\omega^\sharp(\mu)}\omega  \\
        &= \mathcal{L}_{\omega^\sharp(\mu)}(F) \nu + F\mathcal{L}_{\omega^\sharp(\mu)}(\nu) -  F\theta(\omega^\sharp(\mu))\nu - F\mathcal{L}_{\omega^\sharp(\nu)}\mu - dF\wedge \iota_{\omega^\sharp(\nu)}\mu + F\theta(\omega^\sharp(\nu))\mu  \\
        &\qquad+ dF\wedge\iota_{\omega^\sharp(\nu)}\iota_{\omega^\sharp(\mu)}\omega + Fd\iota_{\omega^\sharp(\nu)}\iota_{\omega^\sharp(\mu)}\omega - F\theta\iota_{\omega^\sharp(\nu)}\iota_{\omega^\sharp(\mu)}\omega.
    \end{split}
\end{equation}
Note that the last equality comes from the linearity of $\theta,\omega^\sharp$, and the relation
\begin{equation}
    \mathcal{L}_{F\omega^\sharp(\nu)}\mu = F\mathcal{L}_{\omega^\sharp(\nu)}\mu - dF \wedge \iota_{\omega^\sharp(\nu)}\mu.
\end{equation}
Furthermore, the definition of $\omega^{\flat}$ implies that $\iota_{\omega^\sharp(\mu)} = \mu$. Therefore, we conclude that
\begin{equation}
    \{\mu,F\nu\}_{\Gamma^1(N)} = F\{\mu,\nu\}_{\Gamma^1(N)} + \mathcal{L}_{\omega^\sharp(\mu)}(F) \nu.
\end{equation}

\textbf{Subalgebra of LdR Closed Forms.} A subbundle of the Lie algebroid \eqref{geomdiagram-} can be obtained as follows. At first, according to Cartan's formula, we obtain that
\begin{equation}
\mathcal{L}_{\omega^\sharp(\mu)}\nu -  \theta(\omega^\sharp(\mu))\nu = \iota_{\omega^\sharp(\mu)}d_\theta\nu + d_\theta\iota_{\omega^\sharp(\mu)}\nu , \qquad \mathcal{L}_{\omega^\sharp(\nu)}\mu -  \theta(\omega^\sharp(\nu))\mu = \iota_{\omega^\sharp(\nu)}d_\theta\mu + d_\theta\iota_{\omega^\sharp(\nu)}\mu,
\end{equation}
see also \cite{ChanMurp19}. So, if we suppose that both $\mu$ and $\nu$ are LdR closed one-forms, equation \eqref{Lie-LCS} reduces to
\begin{equation}\label{Lie-LCS2}
\{\mu,\nu\}_{\Gamma^1(N)} = d_\theta \big( \iota_{\omega^\sharp(\mu)}\nu - \iota_{\omega^\sharp(\nu)}\mu + \iota_{\omega^\sharp(\nu)}\iota_{\omega^\sharp(\mu)}\omega \big).
\end{equation}
In particular, let $\mu = d_\theta F$ and $\nu = d_\theta H$, where $F$ and $H$ are real-valued functions defined on $N$. Directly using equation \eqref{Lie-LCS2} and the Hamilton equation \eqref{semiglobal} for $F$ and $H$, we observe that
\begin{equation}
\begin{split}
\{d_\theta F,d_\theta H\}_{\Gamma^1(N)} &= d_\theta \big( \iota_{X_F}d_\theta H - \iota_{X_H}d_\theta F + \iota_{X_H}\iota_{X_F}\omega \big) \\
&= d_\theta \big(\iota_{X_F}\iota_{X_H}\omega - \iota_{X_H}\iota_{X_F}\omega + \iota_{X_H}\iota_{X_F}\omega \big) \\
&= d_\theta \{F,H\}_{LCS}.
\end{split}
\end{equation}
In other words, the set of all LdR exact one-forms on $N$ is a Lie subalgebra of the Lie algebra of one-forms on $N$. This also gives that 
the LdR differential is a Lie algebra homomorphism.  Further, by determining a subbundle of $T^*N$, LdR exact one-forms constitute a Lie subalgebroid of $(T^*N,\pi_N,N,\omega^\sharp,\{\bullet,\bullet\}_{\Gamma^1(N)})$.

\section{Locally Conformal Cosymplectic Manifolds}
\label{Sec-Cos}

 \subsection{Cosymplectic Hamiltonian Dynamics}\label{Sec-Cos-Ham}

An almost cosymplectic structure on a manifold $M$ of odd dimension $2n + 1$ consists of a pair of differential forms $(\eta,\Omega)$, where $\eta$ is a one-form, the so-called  Reeb-form and $\Omega$ is a two-form, and both satisfying the
non-degeneracy condition $\eta\wedge\Omega^n\neq 0$ \cite{Lichnerowicz-Poi,Lichnerowicz-Jacobi}. In other words, the top-form $\eta\wedge\Omega^n$ is a non-vanishing $2n + 1$ volume form on $M$. We refer to a recent and very comprehensive survey \cite{cape} for some history of the theory and for the results we list in the sequel. 

For an almost cosymplectic manifold $(M,\eta,\Omega)$, the mapping
\begin{equation}\label{flat-map}
    \flat: \mathfrak{X}(M)\longrightarrow \Gamma^1(M),\qquad X\mapsto \flat(X)=\iota_X\Omega + \eta(X)\eta
\end{equation}
is an isomorphism. 
We denote the inverse of this isomorphism by $\sharp$. We call the image $\mathcal{R}:=\sharp(\eta)$ the Reeb field. The Reeb vector field fulfills the two following identities that determine the Reeb vector field uniquely.
\begin{equation}\label{Reeb-cosymp}
\iota_{\mathcal{R}}\eta=1,\qquad \iota_{\mathcal{R}}\Omega=0.
\end{equation}

An almost cosymplectic manifold $(M,\eta,\Omega)$ 
turns out to be a cosymplectic manifold if both of the differential forms are closed, i.e., $d\Omega=0$ and $d\eta=0$.

\textbf{Gradient, Hamiltonian and Evolution Vector Fields.}
Consider a Hamiltonian function $H$ defined on an almost cosymplectic manifold $(M,\eta,\Omega)$. We define three vector fields related with the Hamiltonian function. The gradient vector field is defined to be 
\begin{equation}
\operatorname{grad} H:=\sharp(dH),
\end{equation}
see, for example, \cite{Leon-GradOnCosymp}. Notice that, alternatively, we can define gradient vector field as 
\begin{equation}\label{Cos-grad-H}
\iota_{\operatorname{grad} H}\eta=\mathcal{R}(H), \qquad 
\iota_{\operatorname{grad} H}\Omega =dH- \mathcal{R}(H)\eta.
\end{equation}
 For a Hamiltonian function $H$, we define Hamiltonian vector field $X_H$ as 
\begin{equation}\label{Cos-X-H}
\iota_{X_H}\eta=0, \qquad 
\iota_{X_H}\Omega =dH- \mathcal{R}(H)\eta.
\end{equation}
It is immediate to see that the Hamiltonian function $H$ is conserved under the flow generated by the Hamiltonian vector field $X_H$, that is $X_H(H)=0$. One can also define the evolution vector field corresponding to the Hamiltonian function $H$ as the vector field that fulfills the two following identities:
\begin{equation}\label{evolution-H}
	\iota_{E_H}\eta=1, \qquad 
	\iota_{E_H}\Omega =dH- \mathcal{R}(H)\eta.
\end{equation}

\textbf{(Almost) Poisson Bracket.}
Referring to the gradient vector field definition and the two-form $\Omega$, we define a bracket of functions on $M$ as  follows
 \begin{equation}\label{Poisson-bra-coi}
 \{F,H\}_{Cos}:=\Omega(\operatorname{grad} F,\operatorname{grad} H).
 \end{equation}
 It is immediate to observe that we can replace the gradient vector fields in \eqref{Poisson-bra-coi} by Hamiltonian vector fields, so we can alternatively rewrite the bracket as
  \begin{equation}\label{Poisson-bra-coi-1}
 \{F,H\} _{Cos}=\Omega(X_F,X_H)=X_H(F). 
 \end{equation}
One can check that the bracket is an almost Poisson bracket by being skew-symmetric and satisfying the Leibniz identity. 
The bracket \eqref{Poisson-bra-coi} satisfies the Jacobi identity if  and only if $d\Omega=0$ and $d\eta=0$. Then, we deduce that the bracket \eqref{Poisson-bra-coi} is Poisson if and only if the triple $(M,\Omega,\eta)$ is a cosymplectic manifold. 

The following theorem establishes the \textit{symplectization} of a cosymplectic manifold. 
For proof of this assertion, we  cite \cite{Leon-CosympReduction}. 
\begin{proposition}\label{Prop-CS-S}
	Let $M$ be a $(2n+1)$-dimensional manifold equipped with a one-form $\eta$ and a two-form $\Omega$. Let $\overline{M} = M \times \mathbb{R}$ and define the two-form 
\begin{equation}	\label{Omega-bar}
	\overline{\Omega} = \pi^*\Omega+\pi^*\eta \wedge ds,
\end{equation}	
	 where $\pi:\overline{M}\to M $ is the canonical projection and $s$ is the coordinate function on $\mathbb{R}$. Then
	\begin{itemize}
	\item[(i)] $(M,\eta,\Omega)$ is a cosymplectic manifold if and only if $(\overline{M},\overline{\Omega})$ is a symplectic manifold.
	\item[(ii)] $\pi$ is a Poisson morphism.
	\end{itemize}
\end{proposition}

\textbf{Extended Cotangent Bundle.} Generic examples of cosymplectic manifolds are extended cotangent bundles $T^*Q\times \mathbb{R}$. We define a two-form $\Omega$ by pulling back the canonically symplectic two-form $\Omega_Q$ on $T^*Q$ and a one-form $\eta$ pulling back the differential form $dt$ on $\mathbb{R}$. In terms of Darboux coordinates  $(q^i,p_i,t)$ on the bundle, one has
\begin{equation}\label{canonical-cosym}
\Omega=dq^i\wedge dp_i,\qquad \eta=dt.
\end{equation}
A direct computation reads that the triple $(T^*Q\times \mathbb{R},\eta,\Omega)$ turns out to be a cosymplectic manifold.  In this local realization, we compute the Reeb field, the Hamiltonian  and the evolution vector fields as
\begin{equation}\label{loc-X-H-E-H}
\mathcal{R}=\frac{\partial}{\partial t},\qquad 
X_H(q,p,t)=\frac{\partial H}{\partial p_i} 
\frac{\partial }{\partial q^i} 
-
\frac{\partial H}{\partial q^i} 
\frac{\partial }{\partial p_i},\qquad E_H(q,p,t)=\frac{\partial H}{\partial p_i} 
\frac{\partial }{\partial q^i} 
-
\frac{\partial H}{\partial q^i} 
\frac{\partial }{\partial p_i}+\frac{\partial}{\partial t},
\end{equation}
respectively. 
The Poisson bracket in \eqref{Poisson-bra-coi} is computed to be
\begin{equation}
\{F,H\}(q,p,t)=\frac{\partial H}{\partial p_i} 
\frac{\partial F}{\partial q^i} 
-
\frac{\partial H}{\partial q^i} 
\frac{\partial F}{\partial p_i}.
\end{equation} 
and the dynamic generated by the evolution vector field $E_H$ is precisely  
\begin{equation}\label{hamileq22}
 {\dot q}^i =\frac{\partial H}{\partial p_i}, \qquad
 {\dot p}_i =-\frac{\partial H}{\partial q^i}, \qquad
{\dot t} =1.
 \end{equation} 
Since $\dot{t}=1$, we can consider $t$ as a time-parameter up to an affine term. So, we can argue that cosymplectic Hamiltonian dynamics is the convenient geometric formulation for time-dependent Hamiltonians.

Let us present the symplectization of a cosymplectic manifold $( T^*Q\times\mathbb{R})\times \mathbb{R}$ by using Theorem \ref{Prop-CS-S}. The first step here is to determine an extension of this manifold by $\mathbb{R}$, that is $( T^*Q\times\mathbb{R})\times \mathbb{R}$. This is isomorphic to the cotangent bundle $ T^*(Q\times\mathbb{R})\cong T^*Q\times T^* \mathbb{R}$. 
Accordingly, for $z$ in $T^*Q$, we define the following fibrations
\begin{equation}\label{pi-tau}
\begin{split}
\pi&: T^*Q\times T^* \mathbb{R} \longrightarrow T^*Q\times\mathbb{R}, \qquad (z, t, s) \mapsto (z,t),
\\
\tau&: T^*Q\times\mathbb{R} \longrightarrow  Q\times\mathbb{R}, \qquad (z, t) \mapsto (\pi_Q(z),t ),
\end{split}
\end{equation}
Referring to the definition in \eqref{Omega-bar}, in Darboux coordinates, we compute the symplectic two form $\overline{\Omega}$ on $T^*(Q\times\mathbb{R})$ as 
 \begin{equation}\label{Omega-bar-local}
\overline{\Omega}=dq^i\wedge dp_i + dt\wedge ds.
\end{equation}

\textbf{Hamilton-Poincar\'{e} Realization.}
Consider the cotangent bundle $T^*Q$ admitting the Liouville one-form $\theta_Q$ and the symplectic two-form $\Omega_Q$. Recall also the cosymplectic pair $(\eta,\Omega)$ in \eqref{canonical-cosym} defined on the extended cotangent bundle $T^*Q\times \mathbb{R}$.   
An alternative cosymplectic structure on the extended cotangent bundle is determined through the Hamilton-Poincar\'{e} forms
  \begin{equation}\label{thetah}
\theta_H=\Theta_Q-Hdt,\qquad \Omega_H=\Omega_Q+dH\wedge dt,
\end{equation} 
respectively. 
See that, we exhibit the pull back of the canonical form to the extended bundle $T^{*}Q\times \mathbb{R}$ with same notation. 
Recalling $\eta=dt$, a simple observation proves that $(T^{*}Q\times \mathbb{R},\eta,\Omega_H)$ is a cosymplectic manifold. For this second cosymplectic structure the Reeb field $\mathcal{R}_H$ is defined to be
\begin{equation}\label{reebc}
\iota_{\mathcal{R}_H}dt=1,\quad \iota_{\mathcal{R}_H}\Omega_H=0.
\end{equation}
Let us show that the evolution vector field $E_H$ and $\mathcal{R}_H$ are the same. For this, we compare their definitions given in \eqref{evolution-H} and \eqref{reebc}, respectively. The first identities are exactly the same. Let us examine the second identities. We start with the one in \eqref{reebc}, by  substituting the definition of the Hamilton-Pooincar\'{e} two-form in \eqref{thetah}. This gives 
\begin{equation}
\iota_{\mathcal{R}_H}\Omega_Q
=dH-\big(\iota_{\mathcal{R}_H}dH\big)dt
\end{equation}
Recall that, the relationship between the evolution vector field and the Hamiltonian vector field is computed to be $E_H = \mathcal{R} + X_H$. This identity and the conservation of the Hamiltonian function permit us to replace the Reeb field with $E_H$ in the second term of the second identity in \eqref{evolution-H}. So, we argue that $E_H=\mathcal{R}_H$.

\subsection{HJ Theorem for Cosymplectic Hamiltonian Dynamics}\label{Sec-HJ-Cos}

Consider a time-dependent Hamiltonian function $H=H(q,p,t)$ defined on the extended cotangent bundle $T^*Q\times \mathbb{R}$. The Hamilton-Jacobi problem consists on determining a smooth function $S=S(q,t)$ defined on the extended space $Q\times \mathbb{R}$ satisfying 
\begin{equation}\label{HJ-time-dep}
    \frac{\partial S}{\partial t}+H\left(t,q^i,\frac{\partial S(t,q^i)}{\partial t}\right)=f(t).
\end{equation}
Here, $f(t)$ is an arbitrary real-valued real function. In the most classical realization of a time-dependent HJ theory, the function $f$ is taken to be zero.  

Let us comment the ingredients of this theory one by one in order to arrive at the geometric Hamilton -Jacobi theory for the cosymplectic framework. At first, we record here the following sections
\begin{equation}\label{gamma-tilde}
\begin{split}
\Tilde{\gamma}&:Q\times\mathbb{R}\longrightarrow  T^*Q\times T^*\mathbb{R} , \qquad (q,t)\mapsto (\bar{\gamma}(q,t),t,\bar{\bar{\gamma}}(q,t)),
\\
\gamma&:Q\times\mathbb{R}\longrightarrow T^* Q\times\mathbb{R}, \qquad  (q,t)\mapsto (\bar{\gamma}(q,t),t).
\end{split}
\end{equation}
Here, considering $T^*R = \mathbb{R} \times \mathbb{R}$ where $T^*_t\mathbb{R}\cong \mathbb{R}$ for $t \in \mathbb{R}$, $\bar{\gamma}(q,t)$ is in $T^*Q$ whereas,  $\bar{\bar{\gamma}}(q,t)$ is in $\mathbb{R}$. Notice that, in this notation, the following equation holds   
\begin{equation}\label{tildegamma-gamma}
\gamma=\pi\circ \tilde{\gamma}.
\end{equation} 

Consider now a time-dependent Hamiltonian function $H$ on the cosymplectic manifold $T^*Q\times \mathbb{R}$. Notice that, the evolution flow $E_H$ generated by the Hamiltonian function $H$ is defined by the cosymplectic formulation in \eqref{evolution-H} and its local characterization is precisely in the form \eqref{loc-X-H-E-H}. Referring to the projection $\pi$, we pull-back  $H$ and define the extended Hamiltonian 
 \begin{equation}\label{H-s}
H^s(z,t,s)=\pi^*H(z,t)+s
\end{equation}
on the cotangent bundle $T^*(Q\times\mathbb{R})$. 
The flow $X_{H^s}$ generated by the Hamiltonian function $H^s$ is determined through the symplectic relation $\iota_{X_{H^s}}\overline{\Omega} = dH^s$. Here, $\overline{\Omega}$ is the (lifted) symplectic two-form exhibited in \eqref{Omega-bar-local}. In terms of Darboux coordinates, these vector fields are computed to be
\begin{equation} 
X_{H^s}  = \frac{\partial H}{\partial p_i}\frac{\partial}{\partial q^i} -\frac{\partial H}{\partial q^i}\frac{\partial}{\partial p_i} + \frac{\partial}{\partial t} - \frac{\partial H}{\partial t} \frac{\partial}{\partial s}. 
\end{equation}
It is evident from the local observation that $X_{H^s}$ and $E_H$ are $\pi$-related \cite{LeonSar2}, that is,
\begin{equation}\label{pi-related}
E_H \circ \pi = T\pi \circ X_{H^s}.
\end{equation}
We project the vector fields $X_{H^s}$ and $E_H$ to the extended configuration space $Q\times\mathbb{R}$ by means of the sections $\tilde{\gamma}$  and $\gamma$. Accordingly, we arrive at the following vector fields 
\begin{equation}\label{gamma-Ham} 
X^{\Tilde{\gamma}}_{H^s} := T\pi_{Q\times\mathbb{R}} \circ X_{H^s} \circ \Tilde{\gamma},  \qquad 
E^{\gamma}_H  : =  T\tau \circ X_{H^s} \circ \gamma,  
\end{equation}
respectively. The projected vector fields $X^{\Tilde{\gamma}}_{H^s}$ and $E^{\gamma}_H $ coincide. By  choosing $\bar{\gamma}=\bar{\gamma}_i(q,t)dq^i$ we have that 
\begin{equation}\label{gamma-Ham--}
\begin{split}
X^{\Tilde{\gamma}}_{H^s} =E^{\gamma}_H =   \frac{\partial}{\partial t}+  \frac{\partial H}{\partial p_i}\Big\vert_{p_i=\bar{\gamma}_i(q,t)}\frac{\partial}{\partial q^i} . 
\end{split}
\end{equation}
We state the following geometric Hamilton-Jacobi theorem for cosymplectic manifolds \cite{LeonSar2} . 

\begin{theorem}\label{gamma-rel-thm-Cos}
Assume that $\tilde{\gamma}$ is a closed form, then the following conditions are equivalent:
\begin{enumerate}
\item $X^{\Tilde{\gamma}}_{H^s}$ and $E_H$ are $\gamma$-related, that is,
\begin{equation}\label{Cos-gamma-rel}
T\pi \circ T\Tilde{\gamma} \circ X^{\Tilde{\gamma}}_{H^s} = E_H \circ \pi \circ \Tilde{\gamma}.
\end{equation}
\item  The following identiy holds 
\begin{equation}\label{HJ-eq_Cos}
	d(H^s \circ \Tilde{\gamma}) \in \langle \eta \rangle
\end{equation}
where $\langle \eta \rangle $ stands for the span space of the one-form $\eta$. 
\end{enumerate}
\end{theorem}

Let us try to write a geometric Hamilton-Jacobi theory referring to section $\gamma$. For this we need to first introduce some more notation. 
Recall the fibration $\tau$ given in \eqref{pi-tau} and a section $\gamma$ in \eqref{gamma-tilde}. We freeze the base manifold component $q$ and the time variable $t$ one by one to arrive at the following differentiable maps 
\begin{equation} \label{gammas}
\begin{split}
\gamma^{q} &:\mathbb{R}\longrightarrow T^*Q\times \mathbb{R},\qquad t\mapsto \gamma(q,t)=(\bar{\gamma}(q,t),t)
\\
\gamma^t &:Q\longrightarrow T^*Q\times \mathbb{R},\qquad q\mapsto \gamma(q,t)=(\bar{\gamma}(q,t),t)
\end{split}
\end{equation}
where we assume that the following section of the cotangent bundle 
\begin{equation}
\bar{\gamma}^t:Q\longrightarrow T^*Q,\qquad  q\mapsto \bar{\gamma}(q,t)
\end{equation}
is closed for all $t$. Recall the one-form $\tilde{\gamma}$ defined on $Q\times \mathbb{R}$ and examine the implications of its closedness. To have this, we write
\begin{equation*}
\tilde{\gamma}(q,t)=\bar{\gamma}_i (q,t) dq^i + \bar{\bar{\gamma}}(q,t) dt 
\end{equation*}
in terms of the local coordinates. This reads the two next identities:
\begin{equation}\label{local-ident-cosy}
\frac{\partial \bar{\gamma}_i}{\partial q^j}=
\frac{\partial \bar{\gamma}_j}{\partial q^i}, \qquad 
\frac{\partial \bar{\bar{\gamma}}}{\partial q^i}=
\frac{\partial \bar{\gamma}_i}{\partial t}.
\end{equation}
See that the former identity, that expresses the symmetry of the quantity $\bar{\gamma}_{i,j}$ with respect to the change of the indices, manifests that   $\bar{\gamma}^t$ is  a closed one-form. The second is the identity establishing a relationship between $\bar{\gamma}$ and $\bar{\bar{\gamma}}$. 

Let us recall the second condition in Theorem \ref{gamma-rel-thm-Cos}. The left hand side of this is equal to the following expression
\begin{equation*}
\begin{split}
d(H^s\circ \tilde{\gamma})&=d\big(\pi^*H(\bar{\gamma}(q,t),t)+\bar{\bar{\gamma}}(q,t)\big)
\\&=\left(\frac{\partial H}{\partial q^i} 
+\frac{\partial H}{\partial p_j}\frac{\partial \bar{\gamma}_j}{\partial q^i} + \frac{\partial \bar{\bar{\gamma}}}{\partial q^i}
\right) dq^i 
+
\left(\frac{\partial H}{\partial t} + \frac{\partial H}{\partial p_j}\frac{\partial \bar{\gamma}_j}{\partial t} + 
\frac{\partial \bar{\bar{\gamma}}}{\partial t}
\right) dt.
\end{split}
\end{equation*}
The second condition imposes that the one form is in the span of $dt$, reading the identity
\begin{equation*}
\frac{\partial H}{\partial q^i} 
+\frac{\partial H}{\partial p_j}\frac{\partial \bar{\gamma}_j}{\partial q^i} + \frac{\partial \bar{\bar{\gamma}}}{\partial q^i}=0. 
\end{equation*}
See that we can substitute the second identity in \eqref{local-ident-cosy} into the previous equation and arrive at:
\begin{equation*}
\frac{\partial H}{\partial q^i} 
+\frac{\partial H}{\partial p_j}\frac{\partial \bar{\gamma}_j}{\partial q^i}+\frac{\partial \bar{\gamma}_i}{\partial t}=0,
\end{equation*}
that depends only on the Hamiltonian function $H$ and the section $\gamma$. In the literature, this equation is known as the time-dependent Hamilton-Jacobi equation. 

\textbf{Coordinate-free Realization.} 
Assume a differential one-form $\alpha$ on a manifold $Q$. Its vertical lift $\alpha^V$ is a vector field on $T^*Q$ defined to be 
\begin{equation}
\alpha^V: =\omega_Q^\sharp \circ \pi_Q^* \alpha 
\end{equation}
where $\omega_Q$ is the canonical symplectic two-form on $T^*Q$. Here, $\pi_Q^* \alpha$ denotes the pullback of the one-form to the cotangent bundle by means of the cotangent bundle projection.  That is, the following diagram commutes.
\begin{equation}
\xymatrix{ T^*T^{*}Q 
\ar[ddr]^{\pi_{T^*Q}} \ar[rr]^{\omega^\sharp_Q}&     &T T^{*}Q \ar[ddl]_{\tau_{T^*Q}}\\
  &    &\\ &   T^*Q\ar@/^1pc/[uul]^{\pi^*\alpha}
\ar@/_1pc/[uur]_{\alpha^V}}
\end{equation}
We can write this geometrically as follows 
\begin{equation*}
[d(H\circ \gamma^t)]^{V}=\dot{\gamma}^q, 
\end{equation*}
where $\gamma^q$ and $\gamma^t$ are functions in \eqref{gammas}. Here, the superscript $V$ stands for the vertical lift operation resulting with a vector field over the cotangent bundle. 
This permits us to determine an alternative formulation of the Hamilton-Jacobi theorem. In order to see this we recall that the first condition \eqref{Cos-gamma-rel} can also be written as
\begin{equation*}
T\gamma \circ E_H^\gamma = E_H \circ \gamma
\end{equation*}
since $E_H^\gamma$ and $X^{\Tilde{\gamma}}_{H^s} $ coincide.  This reads the following Hamilton-Jacobi theorem.
\begin{theorem}\label{thm-Cristina-cos}
For $\gamma(q,t)=(\bar{\gamma}(q,t),t)$ where $\bar{\gamma}(q,t)$ is closed for each $t$, the following two statements are equivalent
\begin{enumerate}
\item The vector fields $E_H$ and $E_H^{\gamma}$ are $\gamma$-related.
\item The identity $[d(H\circ \gamma^t)]^{V}=\dot{\gamma}^q$ is satisfied.
\end{enumerate}
\end{theorem}
In this theorem, we may further substitute the evolution vector field $E_H$ and its projection $E^{\gamma}_H$ with the Reeb field $\mathcal{R}_H$ and $\mathcal{R}^{\gamma}_H$, respectively.

\subsection{Locally Conformally Cosymplectic Manifolds}

An almost cosymplectic manifold $(M,\eta,\Omega)$ is a locally conformal cosymplectic manifold if there exists a neighborhood $U_\alpha$ of every point $x\in M$ and local a function $\sigma_\alpha: U_\alpha\to \mathbb{R}$ such that
\begin{equation}\label{local-forms-cosymplectic}
d(e^{-\sigma_\alpha}\eta\vert_\alpha) = 0, \qquad
d(e^{-2\sigma_\alpha}\Omega\vert_\alpha) = 0, 
\end{equation}
where $\eta\vert_\alpha$ and $\Omega\vert_\alpha$ denote the restrictions of $\eta$ and $\Omega$ to $U_\alpha$, respectively. See, for example, \cite{LeonLCCos,alpha-cosymp, Leon-CosympReduction}.  According to these definitions, we define a local one-form and a local two-form as
\begin{equation}\label{local-forms}
\eta_\alpha:= e^{-\sigma_\alpha}\eta\vert_\alpha, \qquad \Omega_\alpha := e^{-2\sigma_\alpha}\Omega\vert_\alpha,
\end{equation}
respectively. Evidently, $\eta_\alpha$ and $\Omega_\alpha$ are closed, and they satisfy the non-degeneracy condition $\eta_\alpha \wedge \Omega_\alpha^n\neq 0$. Hence, $(U_\alpha,\eta_\alpha,\Omega_\alpha)$ turns out to be a cosymplectic space.

Now, we consider another local open chart $U_\beta$ such that the intersection $U_\alpha \cap U_\beta$ is not empty. So that, in the intersection, the restriction of the global forms are equal, i.e., $\eta\vert_\alpha = \eta\vert_\beta$ and $\Omega\vert_\alpha = \Omega\vert_\beta$. These identifications and the definitions in \eqref{local-forms} lead us to that, in the overlapping, 
the local one-form and the local two-form are related by  
\begin{equation}\label{alpha-change}
\eta_\beta = e^{-(\sigma_\beta-\sigma_\alpha)}\eta_\alpha, \qquad \Omega_\beta = e^{-2(\sigma_\beta-\sigma_\alpha)}\Omega_\alpha,
\end{equation}
respectively. The relationship between the local one-forms yields a set of scalars
\begin{equation}\label{cocyc-kappa}
\kappa_{\beta\alpha} = e^{-(\sigma_\beta-\sigma_\alpha)} 
\end{equation}
satisfying the cocycle condition in the overlapping of three open neighborhoods. The cocycle character of the local one-forms enables us to glue them up to a line bundle $L\mapsto M$ valued one-form. Similarly, in the overlapping, the local two-forms provide a set of scalars 
\begin{equation}
 \lambda_{\beta\alpha} = e^{-2(\sigma_\beta-\sigma_\alpha)}
 \end{equation}
satisfying the cocycle condition as well. So that, 
the local two-forms $\Omega_\alpha$ can be glued to a line bundle $L\mapsto M$ valued two-form. 

\textbf{The Reeb Field.} 
The non-degeneracy of the pair $(\eta_\alpha,\Omega_\alpha)$ enables us to write the local musical isomorphism
\begin{equation}
\flat_\alpha: \mathfrak{X}(U_\alpha) \to \Gamma(U_\alpha): X_\alpha \mapsto \iota_{X_\alpha}\Omega_\alpha + \eta_\alpha(X_\alpha)\eta_\alpha
\end{equation}
with inverse map $\sharp_\alpha:  \Gamma(U_\alpha) \to \mathfrak{X}(U_\alpha)$. Accordingly, we can write every local one-form $\mu_\alpha$ in $\Gamma(U_\alpha)$ as the image of a unique local vector field $X_\alpha$ under $\flat_\alpha$. 
Referring to the inverse map $\sharp_\alpha$ of the musical isomorphism $\flat_\alpha$, we are able to define the Reeb field $\mathcal{R}_\alpha$ of $U_\alpha$ by $ 
	\mathcal{R}_\alpha = \sharp_\alpha(\eta_\alpha)$.
This definition implies the following local characterizations
\begin{equation}\label{loc-reeb-id}
\iota_{\mathcal{R}_\alpha}\eta_\alpha = 1, \quad \iota_{\mathcal{R}_\alpha}\Omega_\alpha = 0.
\end{equation}
 
Referring to \eqref{local-forms}, we obtain the following local vector fields
\begin{equation}\label{globalReeb-LCCos}
	\mathcal{R}\vert_\alpha = e^{-\sigma_\alpha}\mathcal{R}_\alpha,
\end{equation} 
and a global vector field $\mathcal{R}$ by gluing up all local vector fields $\mathcal{R}\vert_\alpha$. On $U_\alpha \cap U_\beta$, the Reeb fields of $U_\alpha$ and $U_\beta$ are related by $\mathcal{R}_\beta = \kappa_{\alpha\beta}\mathcal{R}_\alpha$, where $\kappa_{\alpha\beta}$ are the scalars given in \eqref{cocyc-kappa}. Therefore, the local Reeb fields $\mathcal{R}_\alpha$ can be glued to a line bundle valued vector field.

\textbf{From Local to Global Picture.} We define the local picture of a global one form and a global vector field as 
\begin{equation}
\mu\vert_\alpha = e^{\sigma_\alpha}\mu_\alpha, \qquad X\vert_\alpha = e^{-\sigma_\alpha} X_\alpha,
\end{equation}
respectively. We compute
\begin{equation*}
\begin{split}
\mu\vert_\alpha &= e^{\sigma_\alpha}\mu_\alpha = e^{\sigma_\alpha}\big(\iota_{X_\alpha}\Omega_\alpha + \eta_\alpha(X_\alpha)\eta_\alpha\big) \\
&= e^{\sigma_\alpha}\big(\iota_{e^{\sigma_\alpha}X\vert_\alpha}e^{-2\sigma_\alpha}\Omega\vert_\alpha + e^{-\sigma_\alpha}\eta\vert_\alpha(e^{\sigma_\alpha}X\vert_\alpha)e^{-\sigma_\alpha}\eta\vert_\alpha\big) \\
& = \iota_{X\vert_\alpha}\Omega\vert_\alpha + \eta\vert_\alpha(X\vert_\alpha)\eta\vert_\alpha.
\end{split}
\end{equation*}
All terms in the last line of this computation can be glued up to their global realizations. Therefore, the gluing up of the local musical isomorphisms $\flat_\alpha$ is precisely the same as the one defined for the manifold $(M,\eta,\Omega)$ with an almost cosymplectic structure.
 
Ignoring the LCC structure for a moment, we can take $(M,\eta,\Omega)$ as an almost cosymplectic manifold. In this setting, we can define the Reeb field of $(M,\eta,\Omega)$ using the identities in \eqref{Reeb-cosymp}. If we take an arbitrary neighborhood $U_\alpha$ and the local Reeb field defined in \eqref{globalReeb-LCCos}, we compute that
\begin{equation*}
\begin{split}
\iota_{\mathcal{R}\vert_\alpha} \eta\vert_\alpha &= \iota_{e^{\sigma_\alpha}\mathcal{R}\vert_\alpha} e^{-\sigma_\alpha}\eta\vert_\alpha = \iota_{\mathcal{R}_\alpha}\eta_\alpha = 1, \\
\iota_{\mathcal{R}\vert_\alpha} \Omega\vert_\alpha &= e^{\sigma_\alpha} \iota_{e^{\sigma_\alpha}\mathcal{R}\vert_\alpha} e^{-2\sigma_\alpha}\Omega\vert_\alpha = 0,
\end{split}
\end{equation*}
that is, the glued Reeb field $\mathcal{R}$ coincides with the Reeb field of an almost cosymplectic manifold $(M,\eta,\Omega)$. 

\textbf{Lee-form Realization.} We can glue up the local definitions in \eqref{local-forms-cosymplectic} in order to arrive at global tensorial formulations. Notice that $d\sigma_\alpha$ can be the local form of a global tensor field $d\sigma_\alpha=d\sigma_\beta$. These locally exact one-forms merge to a closed one-form $\Theta$ on the manifold $M$. The closed one-form $\Theta$ is called the Lee-form. See that the double of these local forms $2d\sigma_\alpha$ determine a global closed form which is twice the Lee-form $2\Theta$. 
Accordingly, we can alternatively say that an almost cosysmplectic manifold $(M,\eta,\Omega)$ is a LCC manifold if
\begin{equation}\label{LCC-str}
d\eta=\Theta\wedge\eta,\qquad d\Omega=2\Theta\wedge\Omega
\end{equation}
for a closed one-form $\Theta$.  In terms of the LdR differential, an almost cosymplectic manifold $(M,\eta,\Omega,\Theta)$ is LCC if and only if the one-form $\eta$ and the two-form $\Omega$ satisfy
\begin{equation}
d_{\Theta}\eta = 0,\qquad d_{2\Theta}\Omega = 0.
\end{equation} 
In this realization, we denote a LCC manifold by a quadruple $(M,\eta,\Omega,\Theta)$ to exhibit the Lee-form $\Theta$. If $\Theta$ is an exact one-form, then we say that $(M,\eta,\Omega,\Theta=d\sigma)$ is a globally conformally cosymplectic (abbreviated as GCCos) manifold. It is evident that if $\Theta=0$, then such a LCC manifold reduces to a cosymplectic manifold.  

\textbf{Extended Cotangent Bundle as a LCC Manifold.}  
Consider the cotangent bundle $T^*Q$ admitting the Liouville one-form $\theta_Q$ and the symplectic two-form $\omega_Q$. We denote the pull-backs of these canonical forms to the extended cotangent bundle $T^*Q\times \mathbb{R}$ by $\Theta_Q$ and $\Omega_Q$, respectively. Assume a closed one-form 
\begin{equation}\label{Lee-Psi}
	\Psi (q,t) =\psi _{j}(q,t)dq^{j}+\zeta(q,t) dt.
\end{equation}
 on the extended configuration space $Q \times \mathbb{R}$ and pull it to back to the extended cotangent bundle $T^*Q\times \mathbb{R}$ by means of the projection $\tau$ in \eqref{pi-tau}. We denote the pull back one-form by $\Theta=\tau^*\Psi$ as well. Define now the following forms
   \begin{equation} \label{LCC-can}
 \eta = d_{\Theta}t=dt-t\Theta  ,\qquad  \Omega_{2\Theta}=-d_{2 \Theta} \Theta _Q=\Omega_Q+2 \Theta\wedge \Theta _Q.
     \end{equation}
It is immediate to check that the quadruple $(T^*Q\times \mathbb{R},\eta,\Omega,\Theta)$ is a LCC manifold.  
In terms of the Darboux coordinates $(q^i,p_i,t)$ on $T^*Q\times \mathbb{R}$, the local realization of this LCC is computed as follows. Assume the 
Lee form in the form
\begin{equation}\label{Lee-LCCos}
	\Theta (q,t) =\psi _{j}(q,t)dq^{j}+\zeta(q,t) dt.
\end{equation}
Then, the differential forms are computed to be
   \begin{equation}
    \begin{split}
    \Omega&=-d_{2 \Theta} \Theta _Q=-d\Theta _Q+2 \Theta\wedge \Theta _Q=dq^i \wedge dp_i+2p_i\Theta \wedge dq^i
    \\
    \eta&= d_{\Theta}t=dt-t\Theta.
            \end{split}
     \end{equation}
We present the following summarizing table for LCC manifolds.
\begin{center}
\begin{tabular}{ ||c||c|c||c||  }
 \hline
 \multicolumn{4}{||c||}{\textbf{LCC Manifolds}} \\
 \hline \hline
 Tensor Field &  Global & Local & Relation\\
 \hline One-form& 
 $\eta\vert_\alpha$ & $\eta_\alpha$ & $\eta_\alpha = e^{-\sigma_\alpha}\eta\vert_\alpha$  \\ Two-form& 
 $\Omega\vert_\alpha$ & $\Omega_\alpha$ & $\Omega_\alpha = e^{-2\sigma_\alpha}\Omega\vert_\alpha$  \\ Ham func & 
 $H\vert_\alpha$ & $H_\alpha$ & $H_\alpha = e^{-\sigma_\alpha}H\vert_\alpha$\\ Ham v-field & 
$X_H\vert_\alpha$ & $X_{H_\alpha}$ & $X_{H_\alpha} = e^{\sigma_\alpha}X_H\vert_\alpha$\\
Evo v-field & 
$E_H\vert_\alpha$ & $E_{H_\alpha}$ & $E_{H_\alpha} = e^{\sigma_\alpha}E_H\vert_\alpha$\\
 \hline Ham Eq & 
 $\iota_{X_H} \eta = 0$ & $\iota_{X_{H_\alpha}} \eta_\alpha = 0$ & \\ &
 $\iota_{X_H} \Omega = d_\Theta H - \langle d_\Theta H,\mathcal{R} \rangle\eta$ & $\iota_{X_{H_\alpha}} \Omega_\alpha = dH_\alpha - \mathcal{R}_\alpha(H_\alpha)\eta_\alpha$ & $\Theta\vert_\alpha = d\sigma_\alpha$\\
 \hline
 Evo Eq & 
 $\iota_{E_H} \eta = 1$ & $\iota_{E_{H_\alpha}} \eta_\alpha = 1$ & \\ &
 $\iota_{E_H} \Omega = d_\Theta H - \langle d_\Theta H,\mathcal{R} \rangle\eta$ & $\iota_{E_{H_\alpha}} \Omega_\alpha = dH_\alpha - \mathcal{R}_\alpha(H_\alpha)\eta_\alpha$ & $\Theta\vert_\alpha = d\sigma_\alpha$\\
 \hline
\end{tabular}
\end{center}
Here, the tensor fields in the global column are those can be glued up to classical (real valued) tensor fields whereas in the local column are those can be glued up to line bundle valued  tensor fields.

\subsection{Relationship Between LCS vs. LCC Manifolds} 	\label{Sec-LCS-LCCos}

Let $(M,\eta,\Omega;\Theta)$ be a $(2n+1)$-dimensional LCC manifold equipped with a Lee-form $\Theta$, and then consider a chart $U_\alpha$ with induced cosymplectic structure $(U_\alpha,\eta_\alpha,\Omega_\alpha)$. Now we apply Proposition \ref{Prop-CS-S} to the local cosymplectic space $(U_\alpha,\eta_\alpha,\Omega_\alpha)$ and arrive at a local symplectic manifold $(U_\alpha\times \mathbb{R},\overline{\Omega}_\alpha)$ where the symplectic two-form is determined to be
\begin{equation}
\overline{\Omega}_\alpha=\pi^*\Omega_\alpha+\pi^*\eta_\alpha\wedge ds_\alpha
\end{equation}
where $s_\alpha$ is the coordinate on the real line. Considering another patch $U_\beta$ and referring to \eqref{alpha-change}, in this chart, we have that
 \begin{equation}
 \begin{split}
\overline{\Omega}_\beta&=\pi^*\Omega_\beta+\pi^*\eta_\beta\wedge ds_\beta \\ &=\pi^*\lambda_{\beta\alpha}\Omega_\alpha+\pi^*\kappa_{\beta\alpha}\eta_\alpha\wedge\pi^* \kappa_{\beta\alpha}ds_\alpha\\ &=\pi^*\lambda_{\beta\alpha}\big( \pi^*\Omega_\alpha+ \pi^*(\eta_\alpha\wedge ds_\alpha) \big)=\pi^*\lambda_{\beta\alpha} \overline{\Omega}_\alpha,
\end{split}
\end{equation}

where we have assumed the local identification 
 \begin{equation}\label{ds-local}
ds_\beta=\pi^*\kappa_{\beta\alpha}ds_\alpha.
\end{equation}
 Let us glue the local extended charts $U_\alpha\times \mathbb{R}$ and obtain the following  manifold
 \begin{equation}\label{M-bar}
\overline{M}:= \bigsqcup_\alpha \big( U_\alpha\times \mathbb{R}\big).
 \end{equation}
 See that, $\overline{M}$ is a one-dimensional fiber bundle over the LCC manifold $M$ and it is not necessarily trivial. We denote this by  $\pi:\overline{M}\mapsto M$. Being the collection of local symplectic spaces, we first examine that whether or not the local symplectic forms can be glued up to a global one on $\overline{M}$. Instead, we notice that, the local character $\overline{\Omega}_\beta=\pi^*\lambda_{\beta\alpha} \overline{\Omega}_\alpha$  provides a global line bundle $\overline{L}\mapsto \overline{M}$ valued two-form  on $\overline{M}$. In order to arrive at global two-forms on $\overline{M}$, we introduce the local form
 \begin{equation}
\overline{\Omega}\big\vert_\alpha:= e^{2(\sigma_\alpha\circ \pi)}\overline{\Omega}_\alpha.
   \end{equation}
A direct observation gives that these local two-forms coincide in the overlapping spaces, that is, $\overline{\Omega}\big\vert_\alpha=\overline{\Omega}\big\vert_\beta$. So that they can be glued up to a well-defined global real valued two-form $\overline{\Omega}$ on $\overline{M}$. In the global view, we have that
 \begin{equation}\label{LCS-Cos}
 d\overline{\Omega}=2\pi^*\Theta \wedge \overline{\Omega}, \qquad  d_{2\pi^*\Theta}\overline{\Omega}=0
    \end{equation}
This reads that $(\overline{M},\overline{\Omega})$ is a locally conformally symplectic manifold with Lee-form $2\pi^*\Theta$. Here, $d_{2\pi^*\Theta}$ stands for the LdR differential. 

See that, the one-forms defined in \eqref{ds-local} can be glued up only to a global line bundle $\overline{L}\mapsto \overline{M}$ valued one-form. On the other hand, we can define a family of local forms
  \begin{equation}\label{u-local}
  u\big\vert_\alpha=e^{\sigma_\alpha\circ \pi}ds_\alpha
      \end{equation}
which, in the light of \eqref{ds-local}, can be glued up to the following one-form 
        \begin{equation}
        u=d_{\pi^*\Theta}s=ds-s\pi^*\Theta.
             \end{equation}
See that this one-form is not closed but it is closed under the LdR differential, that is,
  \begin{equation}
  du=\pi^*\Theta \wedge u, \qquad d_{\pi^*\Theta}u=0.
        \end{equation}
Referring to all these global formulations, we state the following proposition. 
\begin{proposition}\label{LCCos-to-LCS}
 	Consider a LCC manifold $(M,\eta,\Omega)$ with a Lee-form $\Theta$, and then assume a one-dimensional fibration $\pi:\overline{M}\to M$ where $\overline{M}$ is defined as in \eqref{M-bar}. Define a two-form on $\overline{M}$ as
	  \begin{equation}\label{Overline-Omega}
	\overline{\Omega}=\pi^*\Omega+\pi^*\eta\wedge u,
	     \end{equation}
	     where $u$ is the one-form in \eqref{u-local}.  Then, $(M,\eta,\Omega)$ is a LCC manifold with a Lee-form $\Theta$ if and only if $(\overline{M},\overline{\Omega})$ is a LCS manifold with Lee-form $2\pi^*\Theta$.
\end{proposition}

\textbf{Proof.} Let us consider the two-form $\overline{\Omega}$ in \eqref{Overline-Omega} and take the exterior derivative, then
	\begin{equation}
	 \begin{split}
	d\overline{\Omega}&=d\big(\pi^*\Omega+\pi^*\eta\wedge u\big)=
\pi^*d\Omega+\pi^*d\eta\wedge u - \pi^*\eta\wedge d u
\\
&=\pi^*(2\Theta\wedge \Omega)+\pi^*(\Theta\wedge\eta)\wedge u
-\pi^*\eta\wedge \pi^*\Theta \wedge u
\\
&=2\pi^*\Theta\wedge \pi^*\Omega+2\pi^*\Theta \wedge \pi^*\eta \wedge u
\\
&=2\pi^*\Theta\wedge\big(\pi^*\Omega+\pi^*\eta \wedge u\big )
\\
&=2\pi^*\Theta\wedge\overline{\Omega},
\end{split}
	 \end{equation}
	 where we have used the LCC structure identities \eqref{LCC-str} in the second line. 
	This shows that the two-form $\overline{\Omega}$ is a LCS two-form with Lee-form $2\pi^*\Theta$. Compare this with the display in \eqref{LCS-Cos} to see that both the local and the global approaches read the same result. $\blacksquare$
	


	
\textbf{Darboux Coordinates for LCC Manifolds.} 
We have exhibited a Darboux chart for LCS manifolds in the realm of cotangent bundle geometry in Subsection \ref{Sec-LCS}. We have also determined a  
Darboux theorem for cosymplectic manifold in the framework of extended cotangent bundle in Subsection \ref{Sec-Cos-Ham}.  Let us combine these observations with Proposition \ref{LCCos-to-LCS} in order to arrive at a Darboux chart for the LCC manifolds. To this end, we start with an extended cotangent bundle $T^*Q\times \mathbb{R}$, and extend this space with the real line to have $T^*(Q \times \mathbb{R})$. The fibration in this realization is given by  
 \begin{equation}
 \pi:T^*(Q \times \mathbb{R})\cong T^*Q \times T^*\mathbb{R}\longrightarrow T^*Q\times \mathbb{R}, \qquad (z,t,s)\mapsto (z,t).
 \end{equation}
 This is precisely the one in \eqref{pi-tau}. 
 See that $T^*(Q \times \mathbb{R})$ is a symplectic manifold. In Darboux coordinates $(q^i,p_i,t,s)$, the canonical one-form and the symplectic two-form on $T^*(Q \times \mathbb{R})$ can be written as 
  \begin{equation}
  \Theta_{Q \times \mathbb{R}}=p_idq^i+\frac{1}{2}(sdt-tds),\qquad \Omega_{Q \times \mathbb{R}}=-d\Theta_{Q \times \mathbb{R}}=dq^i \wedge dp_i + dt\wedge ds.
   \end{equation}
According to discussions in  Section \ref{Sec-LCS}, one may define a LCS structure on it with a semi-basic closed one-form. We define this form as follows. Take a closed form $\Psi$ on the base manifold $Q \times \mathbb{R}$, and pull this form back to $T^*(Q \times \mathbb{R})$ by means of the cotangent bundle projection $\pi_{Q \times \mathbb{R}}$. We denote the pull back one-form by $\Theta=\pi^*_{Q \times \mathbb{R}}\Psi$. To have a LCS two-form, we simply take minus of LdR differential of the canonical one-form $\Theta_{Q \times \mathbb{R}}$ with the Lee form $2\Theta$ that is
    \begin{equation}\label{bar-omega-1}
    \begin{split}
    \overline{\Omega}&=-d_{2 \Theta} \Theta_{Q \times \mathbb{R}}= -d \Theta_{Q \times \mathbb{R}}+ 2\Theta\wedge \Theta_{Q \times \mathbb{R}}
    \\&=    dq^i \wedge dp_i + dt\wedge ds+ 2p_i\Theta\wedge dq^i+s\Theta\wedge dt- t\Theta\wedge ds.
        \end{split}
     \end{equation}
     This is the canonical form of the LCS two-form $ \overline{\Omega}$ in terms of the Darboux coordinates $(q^i,t,p_i,s)$. 
 On the other hand, recall the symplectic two-form in \eqref{Overline-Omega}  and rewrite it as
  \begin{equation}
      \begin{split}
  \overline{\Omega}&=\pi^*\Omega+\pi^*\eta\wedge u= dq^i \wedge dp_i+2p_i\Theta \wedge dq^i+(dt-t\Theta)\wedge(ds-s\Theta)
   \\&=dq^i \wedge dp_i+dt\wedge ds + 2p_i\Theta \wedge dq^i -sdt\wedge \Theta- t \Theta\wedge ds.
              \end{split}
     \end{equation}
  See that this is precisely equal to the one in \eqref{bar-omega-1}. Let us record this observation in the following proposition saying that every LCC manifold locally looks like the extended cotangent bundle $T^*Q\times \mathbb{R}$ with a Lee-form $\Theta$.  
    \begin{proposition}\label{Darboux-LCC}
	A LCC manifold admits Darboux coordinates $(q^i,p_i,t)$ with a Lee-form $\Theta=\Theta(q,t)$ realizing the LCC structure as  
	  \begin{equation}\label{eta-omega-coord}
    \Omega =dq^i \wedge dp_i+2p_i\Theta \wedge dq^i, \qquad  
    \eta = dt-t\Theta. 
     \end{equation} 
     \end{proposition}


  \section{Locally Conformally Cosymplectic Dynamics and HJ Theory} \label{Sec-Dyn-LCC-Ev}
\subsection{Dynamics on LCC  Manifolds}\label{Sec-Dyn-LCCos}

Let $(M,\eta,\Omega,\Theta)$ be a $(2n+1)$-dimensional LCC manifold equipped with a Lee-form $\Theta$, and consider a chart $U_\alpha$ with induced cosymplectic structure that is $(U_\alpha,\eta_\alpha,\Omega_\alpha)$. By a similar argument to the one given in LCS case, we consider a local Hamiltonian function  $H_\alpha$ on $U_\alpha$ determining a line bundle $L\to M$ valued Hamiltonian function. Furthermore, gluing the local functions
\begin{equation*}
	H\vert_\alpha = e^{\sigma_\alpha}H_\alpha
\end{equation*}
defined on $U_\alpha$, we get a real valued Hamiltonian function $H$ on the whole $M$. We now compute the Hamiltonian, evolution and gradient vector fields one by one.

\textbf{Hamiltonian Vector Field for LCC Manifold.} 
Referring to  \eqref{Cos-X-H} for  cosymplectic spaces, given a local Hamiltonian function  $H_\alpha$ on $U_\alpha$, the local Hamiltonian vector field  $X_{H_\alpha}$ is defined to be
\begin{equation}\label{localHam-LCCos}
\iota_{X_{H_\alpha}}\eta_\alpha=0, \qquad 
\iota_{X_{H_\alpha}}\Omega_\alpha =dH_\alpha- \mathcal{R}_\alpha(H_\alpha)\eta_\alpha.
\end{equation}
By a similar argument as the one given in the LCS case, we    glue the vector fields
\begin{equation}\label{globalHam-LCCos}
	X_H\vert_{\alpha} = e^{-\sigma_\alpha} X_{H_\alpha}
\end{equation}
to get the global Hamiltonian vector field $X_H$. Later, we substitute the local identifications of the differential forms \eqref{local-forms}  and the local characterization of the Reeb field \eqref{globalReeb-LCCos} into the Hamilton equations in \eqref{localHam-LCCos}.  The first equation in \eqref{localHam-LCCos} remains the same, whereas the second equation turns out to be
\begin{equation}
	\begin{split}
		e^{-2\sigma_{\alpha}}\iota_{X_{H_\alpha}}\Omega\vert_\alpha &= dH_\alpha- \mathcal{R}_\alpha(H_\alpha)\eta_\alpha \\
		&= -e^{-\sigma_{\alpha}}d\sigma_{\alpha}H\vert_\alpha + e^{-\sigma_{\alpha}}dH\vert_\alpha - \mathcal{R}\vert_\alpha(e^{-\sigma_{\alpha}})H\vert_\alpha\eta\vert_\alpha - e^{-\sigma_{\alpha}}\mathcal{R}\vert_\alpha(H\vert_\alpha)\\
		&= e^{-\sigma_{\alpha}}\big(-H\vert_\alpha d\sigma_{\alpha} + dH\vert_\alpha + \langle d\sigma_{\alpha},\mathcal{R}\vert_\alpha \rangle H\vert_\alpha \eta\vert_\alpha - \mathcal{R}\vert_\alpha(H\vert_\alpha)\eta\vert_\alpha \big).
	\end{split}	
\end{equation}
We glue up this equation and determine the global Hamilton equations for LCC manifold as
\begin{equation}\label{globalHamEq-LCCos}
\iota_{X_H}\eta=0, \qquad 
\iota_{X_H}\Omega = d_\Theta H - \langle d_\Theta H, \mathcal{R}\rangle \eta.
\end{equation}
For this formulation, we also refer to \cite{LeonLCCos}. It is immediate to see that an alternative way to have the Hamiltonian vector field $X_H$ is
\begin{equation}
	X_H = \sharp (d_\Theta H) - \langle d_\Theta H, \mathcal{R}\rangle\mathcal{R}.
\end{equation}
Here, $\sharp$ is the musical isomorphism obtained through the non-degenerate pair $(\Omega,\eta)$. We note that the pairing of the Hamiltonian vector field $X_H$ and the one-form $d_\theta H$ vanishes identically, that is 
\begin{equation}
\langle d_\Theta H, X_H \rangle=0. 
\end{equation}
This is the locally conformal realization of the conservation of energy. 

\textbf{Evolution Vector Field for LCC Manifolds.} 
To have the equation of motion, we define the local evolution vector field $E_\alpha $ as the sum of the Reeb field  $\mathcal{R}_\alpha$ and the Hamiltonian vector field $X_\alpha$. From equations \eqref{globalReeb-LCCos} and \eqref{globalHam-LCCos}, we can determine a local vector field 
\begin{equation}
	E_H\vert_\alpha = e^{-\sigma_\alpha} E_\alpha,
\end{equation}
then by gluing up these vector fields, we define the global evolution vector field $E_H$ as
\begin{equation}\label{evo-LCC}
	\iota_{E_H}\eta=1, \qquad 
	\iota_{E_H}\Omega = d_\Theta H -\langle d_\Theta H, \mathcal{R}\rangle\eta.
\end{equation}

\textbf{Gradient Vector Field for LCC Manifold.} 
On $(U_\alpha,\eta_\alpha,\Omega_\alpha)$, consider the Hamiltonian function $H_\alpha$. Then, $\operatorname{grad} H_\alpha$ satisfies the equations
\begin{equation}\label{localGrad-LCCos}
\iota_{\operatorname{grad} H_\alpha}\eta_\alpha=\mathcal{R}_\alpha(H_\alpha),\qquad \iota_{\operatorname{grad} H_\alpha}\Omega_\alpha=dH_\alpha - \mathcal{R}_\alpha(H_\alpha)\eta_\alpha.
\end{equation}
From the first equation in \eqref{localGrad-LCCos}, we see that
\begin{equation*}
\begin{split}
	e^{-\sigma_\alpha} \iota_{\operatorname{grad} H_\alpha}\eta\vert_\alpha &= e^{\sigma_\alpha}\mathcal{R}\vert_\alpha(e^{-\sigma_\alpha}H\vert_\alpha) \\
	&= e^{\sigma_\alpha}(-e^{-\sigma_\alpha}\langle d\sigma_\alpha,\mathcal{R}\vert_\alpha\rangle H\vert_\alpha + e^{-\sigma_\alpha}\mathcal{R}\vert_\alpha(H\vert_\alpha)) \\
	&= \mathcal{R}\vert_\alpha(H\vert_\alpha) - \langle d\sigma_\alpha,\mathcal{R}\vert_\alpha\rangle H\vert_\alpha.
\end{split}
\end{equation*}
Therefore,
\begin{equation}\label{eq1-in-localGrad}
	\iota_{e^{-\sigma_\alpha}\operatorname{grad} H_\alpha}\eta\vert_\alpha = \mathcal{R}\vert_\alpha(H\vert_\alpha) - \langle d\sigma_\alpha,\mathcal{R}\vert_\alpha\rangle H\vert_\alpha.
\end{equation}
On the other hand,
\begin{equation*}
\begin{split}
e^{-2\sigma_\alpha}\iota_{\operatorname{grad} H_\alpha}\Omega\vert_\alpha &= d(e^{-\sigma_\alpha}H\vert_\alpha) - e^{\sigma_\alpha}\mathcal{R}\vert_\alpha (e^{-\sigma_\alpha}H\vert_\alpha)e^{-\sigma_\alpha}\eta\vert_\alpha \\
&= -e^{-\sigma_\alpha}H\vert_\alpha d\sigma_\alpha + e^{-\sigma_\alpha}dH\vert_\alpha - \big( -e^{-\sigma_\alpha}\langle d\sigma_\alpha,\mathcal{R}\vert_\alpha\rangle H\vert_\alpha + e^{-\sigma_\alpha}\mathcal{R}\vert_\alpha(H\vert_\alpha)\big)\eta\vert_\alpha \\
&= e^{-\sigma_\alpha} \{ dH\vert_\alpha - H\vert_\alpha d\sigma_\alpha + \big(\langle d\sigma_\alpha,\mathcal{R}\vert_\alpha\rangle H\vert_\alpha - \mathcal{R}\vert_\alpha(H\vert_\alpha)\big)\eta\vert_\alpha \},
\end{split}
\end{equation*}
that is, from the second equation in \eqref{localGrad-LCCos}, we obtain
\begin{equation}\label{eq2-in-localGrad}
	\iota_{e^{-\sigma_\alpha}\operatorname{grad} H_\alpha}\Omega\vert_\alpha = dH\vert_\alpha - H\vert_\alpha d\sigma_\alpha - \big(\mathcal{R}\vert_\alpha(H\vert_\alpha) - \langle d\sigma_\alpha,\mathcal{R}\vert_\alpha\rangle H\vert_\alpha\big)\eta\vert_\alpha.
\end{equation}
Considering \eqref{eq1-in-localGrad} and \eqref{eq2-in-localGrad}, we are able to define
\begin{equation}
(\operatorname{grad} H)\vert_\alpha = e^{-\sigma_\alpha}\operatorname{grad} H_\alpha.
\end{equation} 
The global realizations of each term in \eqref{eq1-in-localGrad} and \eqref{eq2-in-localGrad} allow us to write the following equations for $\operatorname{grad} H$ 
\begin{equation}
	\iota_{\operatorname{grad} H}\eta=\langle d_\Theta H, \mathcal{R}\rangle, \qquad 
	\iota_{\operatorname{grad} H}\Omega = d_\Theta H - \langle d_\Theta H, \mathcal{R}\rangle\eta.
\end{equation}
An equivalent definition for the corresponding gradient vector field of a Hamiltonian function $H$ might be given by 
\begin{equation}\label{grad-LCC}
	\operatorname{grad} H = \sharp (d_\Theta H).
\end{equation}
Then, the relation between the gradient vector field and the Hamiltonian vector field becomes
\begin{equation}\label{grad-Ham}
	\operatorname{grad} H - X_H = \langle d_\Theta H, \mathcal{R}\rangle\mathcal{R}.
\end{equation}

Notice that when we focus on the almost cosymplectic structure of $(M,\eta,\Omega)$, $X_H$ in \eqref{globalHamEq-LCCos} does not give the Hamiltonian vector field corresponding to $H$ unless $\Theta = 0$. Accordingly, $E_H$ in \eqref{evo-LCC} (respectively, $\operatorname{grad}H$ in \eqref{grad-LCC}) is not the same as the evolution vector field (respectively, gradient vector field) of the almost cosymplectic manifold $(M,\eta,\Omega)$ unless $\Theta = 0$.

\textbf{Coordinate realizations.} Consider the almost symplectic manifold $(T^*Q\times \mathbb{R},\eta,\Omega,\Theta)$ with Darboux coordinates $(q^i,p_i,t)$.  Assume the Lee form $\theta=\theta(q,t)$ as in \eqref{Lee-LCCos}. The Reeb field $\mathcal{R}$ is computed to be
\begin{equation}\label{local-Reeb-LCC}
\mathcal{R} = \frac{1}{1-t\zeta}\frac{\partial}{\partial t} + \frac{2p_i\zeta}{1-t\zeta}\frac{\partial}{\partial p_i}.
\end{equation}
It is immediate to see that if $\zeta$ (that is the time component of the Lee-form) is zero, then this Reeb field reduces to the classical Reeb field $\partial/\partial t$. So that, in the light of this local calculation we see that the $\zeta$ term in the Lee-form is responsible for the gluing of the time parameter. A direct computation gives the Hamiltonian vector field $X_H$ as
\begin{equation}\label{Ham-LCC-local}
X_H=\frac{\partial H}{\partial p_i} \frac{\partial}{\partial q^i} + \left( -\frac{\partial H}{\partial q^{i}}+\frac{\psi _{i}H}{1-t\zeta }+
\frac{2\left( p_{i}\psi _{j}-p_{j}\psi _{i}\right) }{\left(
1-t\zeta \right) }\frac{\partial H}{\partial p_{j}} - \frac{t\psi _{i}}{1-t\zeta}\frac{\partial H}{\partial t}\right) \frac{\partial}{\partial p_i} + \frac{t\psi _{i}}{1-t\zeta}\frac{\partial H}{\partial p_i} \frac{\partial}{\partial t}.  
\end{equation}

Since the evolution vector field $E_H$ is defined as the sum of the Reeb field $\mathcal{R}$ and the Hamiltonian vector field $X_H$, we can immediately calculate $E_H$ as
\begin{equation}\label{E-H-LCC-local}
\begin{split}
E_H=\frac{\partial H}{\partial p_i} \frac{\partial}{\partial q^i} &+ \left( -\frac{\partial H}{\partial q^{i}}+\frac{\psi _{i}H}{1-t\zeta }+
\frac{2\left( p_{i}\psi _{j}-p_{j}\psi _{i}\right) }{\left(
1-t\zeta \right) }\frac{\partial H}{\partial p_{j}} - \frac{t\psi _{i}}{1-t\zeta}\frac{\partial H}{\partial t} + \frac{2p_i\zeta}{1-t\zeta} \right) \frac{\partial}{\partial p_i} \\
&+ \frac{1}{1-t\zeta}\left(1+t\psi _{i}\frac{\partial H}{\partial p_i}\right) \frac{\partial}{\partial t}. 
\end{split}
\end{equation}
Similarly, the relation in \eqref{grad-Ham} helps us  calculate the gradient vector field $\operatorname{grad}H$ as
\begin{equation}
\begin{split}
\operatorname{grad}H=\frac{\partial H}{\partial p_i} \frac{\partial}{\partial q^i} &+ \left( -\frac{\partial H}{\partial q^{i}}+\frac{\psi _{i}H}{1-t\zeta } - \frac{2p_iH{\zeta}^2}{(1-t\zeta)^2} +  \left(
\frac{2 ( p_{i}\psi _{j}-p_{j}\psi _{i} ) }{ (
1-t\zeta  ) } + \frac{2p_ip_j{\zeta}^2}{(1-t\zeta)^2} \right) \frac{\partial H}{\partial p_{j}} \right. \\
& \left. \qquad+ \left(\frac{2p_i\zeta}{(1-t\zeta)^2} - \frac{t\psi _{i}}{1-t\zeta} \right)\frac{\partial H}{\partial t}\right) \frac{\partial}{\partial p_i} \\
 &+ \left( \left(\frac{2p_i\zeta}{(1-t\zeta)^2} + \frac{t\psi _{i}}{1-t\zeta} \right)\frac{\partial H}{\partial p_i} + \frac{1}{(1-t\zeta)^2}\frac{\partial H}{\partial t} - \frac{H\zeta}{(1-t\zeta)^2}\right) \frac{\partial}{\partial t}. 
\end{split}
\end{equation}
\subsection{Jacobi Structure of LCC Manifolds}\label{sec-Jac-LCC}



Consider a LCC manifold $(M,\eta,\Omega,\Theta)$.   For each local chart $U_\alpha$, the local isomorphism $\sharp_\alpha$ and the local  two-form $\Omega_\alpha$ give a Poisson bivector field $\Lambda_\alpha$ on $U_\alpha$. We write this as  
\begin{equation}\label{bivec-omega}
\Lambda_\alpha(\mu,\nu)=\Omega_\alpha(\sharp_\alpha(\mu),\sharp_\alpha(\nu)). 
\end{equation}
Accordingly, for two local functions $F_\alpha$ and $H_\alpha$, we may define the associated local Poisson bracket by
\begin{equation}\label{bracket-alpha}
\{F_\alpha,H_\alpha\}_\alpha := \Lambda_\alpha(dF_\alpha,dH_\alpha).
\end{equation}
Following the policies determined in the previous section, let us glue up this local realization. For this, we first substitute \eqref{bivec-omega} in \eqref{bracket-alpha} to obtain
\begin{equation} \label{addss}
\begin{split}
\{F_\alpha,H_\alpha\}_\alpha &= \Omega_\alpha(\sharp_\alpha(dF_\alpha),\sharp_\alpha(dH_\alpha)) = \Omega_\alpha(\operatorname{grad}F_\alpha,\operatorname{grad}H_\alpha) = \Omega_\alpha(X_{F_\alpha},X_{H_\alpha}) \\
&= e^{-2\sigma_\alpha} \Omega\vert_\alpha(e^{\sigma_\alpha}X_F\vert_\alpha,e^{\sigma_\alpha}X_H\vert_\alpha) \\
&= \Omega\vert_\alpha(X_F\vert_\alpha,X_H\vert_\alpha).
\end{split}
\end{equation}
where we have used the local identification of the local cosymplectic two-form presented in \eqref{local-forms} and the local identification of the Hamiltonian vector fields in \eqref{globalHam-LCCos}. Notice that the right hand side is global, so is the left hand side. 

\textbf{LCC as a Jacobi Manifold.} 
Let us turn back to the local Poisson bracket \eqref{bracket-alpha} and try to examine it from  its local bivector formulation. The definition in \eqref{bracket-alpha} brings out the direct calculation
\begin{equation}\label{loc-Poiss-global}
\begin{split}
\{F_\alpha,H_\alpha\}_\alpha &= \Lambda_\alpha(dF_\alpha,dH_\alpha) =  \Lambda_\alpha\big(d(e^{-\sigma_\alpha}F\vert_\alpha),d(e^{-\sigma_\alpha}H\vert_\alpha)\big) \\
&= \Lambda_\alpha\big( e^{-\sigma_\alpha}dF\vert_\alpha - e^{-\sigma_\alpha}F\vert_\alpha d\sigma_\alpha, e^{-\sigma_\alpha}dH\vert_\alpha - e^{-\sigma_\alpha}H\vert_\alpha d\sigma_\alpha \big) \\
&= e^{-2\sigma_\alpha}\Lambda_\alpha\big(dF\vert_\alpha - F\vert_\alpha d\sigma_\alpha, dH\vert_\alpha - H\vert_\alpha d\sigma_\alpha \big) \\
&= e^{-2\sigma_\alpha}\Lambda_\alpha\big(dF\vert_\alpha, dH\vert_\alpha \big) - e^{-2\sigma_\alpha}H\vert_\alpha\Lambda_\alpha\big(dF\vert_\alpha,  d\sigma_\alpha \big) - e^{-2\sigma_\alpha}F\vert_\alpha\Lambda_\alpha\big( d\sigma_\alpha, dH\vert_\alpha \big).
\end{split}
\end{equation}
The left hand side is global according to the calculation in \eqref{addss}, this suggests the local identification
\begin{equation}\label{LCC-bi-vector}
\Lambda\vert_\alpha := e^{-2\sigma_\alpha}\Lambda_\alpha.
\end{equation}
We can define a vector field 
\begin{equation}\label{V-theta}
V\vert_\alpha:=\iota_{2d{\sigma_\alpha}}\Lambda\vert_\alpha
\end{equation}
and then state the following proposition.

\begin{proposition}\label{prop-LCC-Jacobi}
A LCC manifold $(M,\eta,\Omega,\Theta)$ admits a Jacobi manifold structure determined by the pair $(V,\Lambda)$ where $V$ is the  vector field in \eqref{V-theta} and $\Lambda$ is the bivector obtained by gluing up the local bivectors as in \eqref{LCC-bi-vector}. 
\end{proposition}

\textbf{Proof.} To show that $(M,V,\Lambda)$ is a Jacobi manifold we need to justify the identities in \eqref{ident-Jac}. 
We start with the local Poisson bivector $\Lambda_\alpha$. Since it is Poisson, its Schouten-Nijenhuis bracket is zero, i.e, $[\Lambda_\alpha,\Lambda_\alpha]=0$. We start with this identity and do the following calculation
\begin{equation}
\begin{split}
[\Lambda_\alpha,\Lambda_\alpha] &= [e^{2\sigma_\alpha}\Lambda\vert_\alpha,e^{2\sigma_\alpha}\Lambda\vert_\alpha] 
\\
&= [e^{2\sigma_\alpha}\Lambda\vert_\alpha,e^{2\sigma_\alpha}]\wedge\Lambda\vert_\alpha + e^{2\sigma_\alpha}[e^{2\sigma_\alpha}\Lambda\vert_\alpha,\Lambda\vert_\alpha] 
\\
&= e^{2\sigma_\alpha}[\Lambda\vert_\alpha,e^{2\sigma_\alpha}] \wedge \Lambda\vert_\alpha + e^{4\sigma_\alpha}[\Lambda\vert_\alpha,\Lambda\vert_\alpha] + e^{2\sigma_\alpha}[e^{2\sigma_\alpha},\Lambda\vert_\alpha] \wedge \Lambda\vert_\alpha 
\\
&= e^{2\sigma_\alpha}(-\iota_{d(e^{2\sigma_\alpha})}\Lambda\vert_\alpha) \wedge \Lambda\vert_\alpha+ e^{4\sigma_\alpha}[\Lambda\vert_\alpha,\Lambda\vert_\alpha] +  e^{2\sigma_\alpha}(-\iota_{d(e^{2\sigma_\alpha})}\Lambda\vert_\alpha)\wedge \Lambda\vert_\alpha 
\\
&= 2e^{4\sigma_\alpha}(-\iota_{2d{\sigma_\alpha}}\Lambda\vert_\alpha)\wedge \Lambda\vert_\alpha + e^{4\sigma_\alpha}[\Lambda\vert_\alpha,\Lambda\vert_\alpha].
 \end{split}
\end{equation}
Since this expression is zero, we can conclude that
\begin{equation}
[\Lambda\vert_\alpha,\Lambda\vert_\alpha] =2 \iota_{2d{\sigma_\alpha}}\Lambda\vert_\alpha\wedge \Lambda\vert_\alpha.
\end{equation}
Notice that all the tensorial fields in this expression admit global realizations. That gives us the following global equality
\begin{equation}\label{LCS-Jac-1}
[\Lambda,\Lambda]=2V\wedge \Lambda
\end{equation}
where we have employed \eqref{V-theta}. This satisfies the first requirement in definition \eqref{ident-Jac}. For the second one, we need to recall some of the properties from \cite{Marle-SN}. First, we compute that
  \begin{equation}\label{Cals-1}
 [V,\Lambda]=[\iota_{2\Theta}\Lambda,\Lambda]=\frac{1}{2} \iota_{2\Theta} [\Lambda,\Lambda],
  \end{equation}
where the second equality is the Schouten-Nijenhuis bracket for any closed one-form $2\Theta$ and any bivector field $\Lambda$. Further, for a closed one-form $2\Theta$, the distribution of the interior product over the Schouten-Nijenhuis bracket reads 
  \begin{equation}\label{Cals-2}
\iota_{2\Theta}[\Lambda,\Lambda]=\iota_{2\Theta} \Lambda \wedge \iota_{2\Theta} \Lambda + \big(\iota_{2\Theta}\iota_{2\Theta}\Lambda\big) \Lambda=0,
  \end{equation}
where, in the last calculation, both terms vanish identically due to the antisymmetric properties. 
Merging the calculations in \eqref{Cals-1} and \eqref{Cals-2}, we conclude that $[V,\Lambda]=0$. So, we can argue now that $(M,V,\Lambda)$ is indeed a Jacobi manifold. $\blacksquare$

We remark here that the previous result provides an example of the so-called locally conformally Poisson manifolds and it fits the geometry presented in \cite{Vais-Dirac,Vita16}.  
It is immediate to apply the Jacobi bracket in Appendix \ref{Sec-Jac-M} to the present case as
\begin{equation}\label{Pois-LCC} 
\{F,H\}_{LCC} = \Lambda(dF,dH) + FV(H) - HV(F). 
\end{equation}

\begin{proposition}\label{Prop-Jacobi-Map}
Consider a LCC  manifold $(M,\eta,\Omega,\Theta)$ and assume that $\pi:\overline{M}\to M$ is a one-dimensional fibration where $(\overline{M},\overline{\Omega};\Theta)$ is the LCS manifold defined as in \eqref{M-bar}, then, $\pi$ is a Jacobi map. 
\end{proposition}

\textbf{Proof.} Notice that in each local chart $U_\alpha$ we have a cosymplectic space $(U_\alpha,\Omega_\alpha,\eta_\alpha)$ and a symplectic space $(U_\alpha \times \mathbb{R},\overline{\Omega}_\alpha)$. In the light of Proposition \ref{Prop-CS-S}, these structures are related with the equation in \eqref{Omega-bar}. Moreover the projection mapping $\pi:U_\alpha \times \mathbb{R}\mapsto U_\alpha$ is a Poisson map. In terms of the local Poisson bivectors this gives that 
\begin{equation}
\pi_*\overline{\Lambda}_\alpha=\Lambda _\alpha.
\end{equation}
Here $\overline{\Lambda}_\alpha$ is the Poisson bivector for the local symplectic structure whereas $\Lambda _\alpha$ is the local Poisson bivector for the local cosymplectic structure. 
We substitute the identification in \eqref{LCC-bi-vector} and obtain the global identification $\pi_*\overline{\Lambda}=\Lambda$ where $\overline{\Lambda}$. Further, from \eqref{V-theta}, one has that 
\begin{equation} 
V\vert_\alpha:=\iota_{2d{\sigma_\alpha}}\Lambda\vert_\alpha = \iota_{2d{\sigma_\alpha}}\pi_*\overline{\Lambda}\vert_\alpha =
\pi_* \iota_{2\pi^*d\sigma_\alpha}\overline{\Lambda}\vert_\alpha =\overline{\Omega}\vert_\alpha ^\sharp (2\pi^*d\sigma_\alpha) = Z_{2\pi^*d\sigma_\alpha }\vert_\alpha.
\end{equation}
$\blacksquare$

\subsection{Algebra of One-forms}\label{Sec-LCC-Algebroid}

Consider a LCC manifold $(M,\eta,\Omega,\Theta)$ and the musical isomorphism $\flat$ defined in \eqref{flat-map} together with its inverse $\sharp$. Let us denote these musical isomorphisms as
\begin{equation}
\flat(X) = X^\flat, \qquad \sharp(\mu) = \mu^\sharp
\end{equation}
for an arbitrary vector field $X$ and an arbitrary one-form $\mu$. By definition of $\flat$, we are able to write the identity
\begin{equation}\label{pair-id2}
\Omega(\mu^\sharp, X) = \langle \mu, X \rangle - \eta(\mu^\sharp)\eta(X).
\end{equation}
We introduce a bracket on the space $\Gamma^1(M)$ of one form sections as follows 
\begin{equation}\label{Lie-bra-LCC}
\{\bullet , \bullet\}_{\Gamma^1(M)}: \Gamma^1(M) \times \Gamma^1(M) \longrightarrow \Gamma^1(M), \qquad (\mu,\nu)\mapsto \flat([\mu^\sharp,\nu^\sharp]),
\end{equation}
where the bracket inside the musical mapping $\flat$ is the Jacobi-Lie bracket of vector fields on $M$. This definition implies that the musical mappings $\flat$ and $\sharp$ are Lie algebra isomorphisms, that is, they satisfy the identities given in \eqref{iso-id}. Mimicking the procedure done for LCS case in Section \ref{Sec-LCS-Algebroid}, we take the Lichnerowicz-deRham differential of $\Omega$ at $(\mu^\sharp,\nu^\sharp,X)$. This gives 
\begin{equation}\label{LDR-omega2}
\begin{split}
d_{2\Theta}\Omega (\mu^\sharp,\nu^\sharp,X) &= d\Omega (\mu^\sharp,\nu^\sharp,X) - 2\Theta\wedge \Omega (\mu^\sharp,\nu^\sharp,X)\\
&= \mu^\sharp\big(\Omega(\nu^\sharp,X)\big) - \nu^\sharp\big(\Omega(\mu^\sharp,X)\big) + X\big(\Omega (\mu^\sharp,\nu^\sharp)\big) \\
& \qquad -\Omega\big([\mu^\sharp,\nu^\sharp],X\big) + \Omega\big([\mu^\sharp,X],\nu^\sharp\big) - \Omega\big([\nu^\sharp,X],\mu^\sharp\big) \\
& \qquad -2\Theta(\mu^\sharp)\Omega(\nu^\sharp,X) + 2\Theta(\nu^\sharp)\Omega(\mu^\sharp,X) - 2\Theta(X)\Omega(\mu^\sharp,\nu^\sharp) \\
&= \mu^\sharp\big( \langle\nu,X\rangle \big) - \mu^\sharp(\eta(\nu^\sharp))\eta(X) - \eta(\nu^\sharp)\mu^\sharp(\eta(X)) \\
& \qquad- \nu^\sharp\big( \langle\mu,X\rangle \big) + \nu^\sharp(\eta(\mu^\sharp))\eta(X) + \eta(\mu^\sharp)\nu^\sharp(\eta(X)) \\
& \qquad+\langle d\iota_{\nu^\sharp}\iota_{\mu^\sharp}\Omega, X\rangle -\langle\{\mu,\nu\}_{\Gamma^1},X\rangle + \eta([\mu^\sharp,\nu^\sharp])\eta(X) \\
& \qquad - \langle \nu,[\mu^\sharp,X]\rangle + \eta(\nu^\sharp)\eta([\mu^\sharp,X]) + \langle \mu,[\nu^\sharp,X]\rangle - \eta(\mu^\sharp)\eta([\nu^\sharp,X])\\
&\qquad -2\Theta(\mu^\sharp)\big(\langle\nu,X\rangle - \eta(\nu^\sharp)\eta(X)\big) + 2\Theta(\nu^\sharp)\big(\langle\mu,X\rangle - \eta(\mu^\sharp)\eta(X)\big) \\
& \qquad - 2\iota_{\nu^\sharp}\iota_{\mu^\sharp}\Omega\langle\Theta,X\rangle.
\end{split}
\end{equation}
Here we have used the identity \eqref{pair-id2}. The commutation property of the Lie derivative and the interior derivative implies that
\begin{equation}
\mu^\sharp\big( \langle\eta,X\rangle \big) -\langle \eta,[\mu^\sharp,X]\rangle = \langle\mathcal{L}_{\mu^\sharp}\eta,X\rangle, \qquad \nu^\sharp\big( \langle\eta,X\rangle \big) -\langle \eta,[\nu^\sharp,X]\rangle = \langle\mathcal{L}_{\nu^\sharp}\eta,X\rangle.
\end{equation}
After we substitute these identities in equation \eqref{LDR-omega}, and since $\Omega$ is LdR closed, that is, $d_{2\Theta}\Omega = 0$, we arrive at the following proposition. 
\begin{proposition}
For a LCC manifold $(M,\eta,\Omega,\Theta)$, the space $\Lambda^1(M)$ of one-form sections is a Lie algebra that admits the Lie bracket 
\begin{equation}\label{Lie-LCC}
\begin{split}
 \{\mu,\nu\}_{\Gamma^1(M)} &= \mathcal{L}_{\mu^\sharp}\nu -  2\Theta(\mu^\sharp)\nu - \mathcal{L}_{\nu^\sharp}\mu + 2\Theta(\nu^\sharp)\mu + d_{2\Theta}\iota_{\nu^\sharp}\iota_{\mu^\sharp}\Omega \\
 &\qquad+ \eta(\mu^\sharp)\big( \mathcal{L}_{\nu^\sharp}\eta - 2\Theta(\nu^\sharp)\eta \big) - \eta(\nu^\sharp)\big( \mathcal{L}_{\mu^\sharp}\eta - 2\Theta(\mu^\sharp)\eta \big) \\
 &\qquad+\big( \eta([\mu^\sharp,\nu^\sharp]) - \mu^\sharp(\eta(\nu^\sharp)) + \nu^\sharp(\eta(\mu^\sharp)) \big)\eta.
\end{split}
\end{equation}
\end{proposition}
For $\Theta=0$, one arrives at the algebra of 
differential one-forms for a cosymplectic manifold. 
In this case, the bracket \eqref{Lie-LCS} reduces to the bracket
 \begin{equation}\label{Lie-LCC-Sym}
\begin{split}
 \{\mu,\nu\}_{\Gamma^1(M)} &=  \mathcal{L}_{\mu^\sharp}\nu  - \mathcal{L}_{\nu^\sharp}\mu  + d \iota_{\nu^\sharp}\iota_{\mu^\sharp}\Omega  + \eta(\mu^\sharp) \mathcal{L}_{\nu^\sharp}\eta    - \eta(\nu^\sharp) \mathcal{L}_{\mu^\sharp}\eta  \\
 & \qquad +\big( \eta([\mu^\sharp,\nu^\sharp]) - \mu^\sharp(\eta(\nu^\sharp)) + \nu^\sharp(\eta(\mu^\sharp)) \big)\eta.
\end{split}
\end{equation}

\textbf{Lie Algebroid Realization.} For a LCC manifold $(M,\eta,\Omega,\Theta)$, we have derived an algebra structure \eqref{Lie-LCC} on the space $\Gamma^1(M)$ of one-form sections. This permits us to realize the Lie algebroid structure $(T^*M,\pi_M,M,\sharp,\{\bullet,\bullet\}_{\Gamma^1(M)})$ by the following  commutative diagram, where the anchor map is the musical isomorphism $\sharp$.
 \begin{equation}\label{geomdiagram----}
\xymatrix{ (T^{*}M,\{\bullet,\bullet\}_{\Gamma^1(M)})
\ar[ddr]_{\pi_M} \ar[rr]^{\quad  \sharp}&   & TM \ar[ddl]^{\tau_M}\\
  &    &\\
 & M}
\end{equation}
 It is straightforward to show that the bracket \eqref{Lie-LCC} satisfies the identity \eqref{oid-id}. Let $F\in \mathcal{F}(M)$ and $\mu,\nu \in \Gamma^1(M)$. We have
 \begin{equation}
     \begin{split}
 \{\mu,F\nu\}_{\Gamma^1(M)} &= \mathcal{L}_{\mu^\sharp}(F\nu) -  2\Theta(\mu^\sharp)F\nu - \mathcal{L}_{(F\nu)^\sharp}\mu + 2\Theta((F\nu)^\sharp)\mu + d_{2\Theta}\iota_{(F\nu)^\sharp}\iota_{\mu^\sharp}\Omega \\
 &\qquad+ \eta(\mu^\sharp)\big( \mathcal{L}_{(F\nu)^\sharp}\eta - 2\Theta((F\nu)^\sharp)\eta \big) - \eta((F\nu)^\sharp)\big( \mathcal{L}_{\mu^\sharp}\eta - 2\Theta(\mu^\sharp)\eta \big) \\
 &\qquad+\big( \eta([\mu^\sharp,(F\nu)^\sharp]) - \mu^\sharp(\eta((F\nu)^\sharp)) + (F\nu)^\sharp(\eta(\mu^\sharp)) \big)\eta \\
 &= \mathcal{L}_{\mu^\sharp}(F)\nu + F\mathcal{L}_{\mu^\sharp}(\nu) -  2F\Theta(\mu^\sharp)\nu - F\mathcal{L}_{\nu^\sharp}\mu - dF \wedge \iota_{\nu^\sharp}\mu + 2F\Theta(\nu^\sharp)\mu \\
 &\qquad+ dF\wedge\iota_{\nu^\sharp}\iota_{\mu^\sharp}\Omega + Fd\iota_{\nu^\sharp}\iota_{\mu^\sharp}\Omega -2F\Theta\wedge \iota_{\nu^\sharp}\iota_{\mu^\sharp}\Omega \\
  &\qquad+ \eta(\mu^\sharp)\big( F\mathcal{L}_{\nu^\sharp}\eta + dF\wedge\iota_{\nu^\sharp}\eta - 2F\Theta(\nu^\sharp)\eta \big) - F\eta(\nu^\sharp)\big( \mathcal{L}_{\mu^\sharp}\eta - 2\Theta(\mu^\sharp)\eta \big) \\
  &\qquad+\big( \mu^\sharp(F)\eta(\nu^\sharp) + F\eta([\mu^\sharp,\nu^\sharp]) - \mu^\sharp(F)\eta(\nu^\sharp) - F\mu^\sharp(\eta(\nu^\sharp)) + F\nu^\sharp(\eta(\mu^\sharp)) \big)\eta \\
  &= \mathcal{L}_{\mu^\sharp}(F)\nu + F\mathcal{L}_{\mu^\sharp}(\nu) -  2F\Theta(\mu^\sharp)\nu - F\mathcal{L}_{\nu^\sharp}\mu + 2F\Theta(\nu^\sharp)\mu + Fd_{2\Theta}\iota_{\nu^\sharp}\iota_{\mu^\sharp}\Omega \\
  &\qquad+ F\eta(\mu^\sharp)\big( \mathcal{L}_{\nu^\sharp}\eta  - 2\Theta(\nu^\sharp)\eta \big) - F\eta(\nu^\sharp)\big( \mathcal{L}_{\mu^\sharp}\eta - 2\Theta(\mu^\sharp)\eta \big) \\
  &\qquad+F\big( \eta([\mu^\sharp,\nu^\sharp])  - \mu^\sharp(\eta(\nu^\sharp)) + \nu^\sharp(\eta(\mu^\sharp)) \big)\eta \\
  & \qquad - dF \wedge \iota_{\nu^\sharp}\mu + dF\wedge\iota_{\nu^\sharp}\iota_{\mu^\sharp}\Omega + \eta(\mu^\sharp)dF\wedge\iota_{\nu^\sharp}\eta
\end{split}
 \end{equation}
 where we have made use of the linearity of $\Theta, \sharp, \eta$, as well as the relations
 \begin{equation}
    \mathcal{L}_{F\nu^\sharp}\mu = F\mathcal{L}_{\nu^\sharp}\mu + dF\wedge\iota_{\nu^\sharp}\mu, \qquad \mathcal{L}_{F\nu^\sharp}\eta = F\mathcal{L}_{\nu^\sharp}\eta + dF\wedge\iota_{\nu^\sharp}\eta, \qquad [\mu^\sharp,F\nu^\sharp] = \mu^\sharp(F)\nu^\sharp + F[\mu^\sharp,\nu^\sharp].
 \end{equation}
 Since the definition of $\flat$ implies that $\iota_{\mu^\sharp}\Omega + \eta(\mu^\sharp)\eta = \mu$, we conclude that
 \begin{equation}
     \{\mu,F\nu\}_{\Gamma^1(M)} = F\{\mu,\nu\}_{\Gamma^1(M)} + \mathcal{L}_{\mu^\sharp}(F)\nu.
 \end{equation}

\subsection{Geometric HJ for LCC Dynamics} \label{Sec-HJ-LCC}

Recall the fibrations and the sections presented in Section \ref{Sec-HJ-Cos}. We are interested again in the extended cotangent bundle $T^*Q\times \mathbb{R}$ and the cotangent bundle $T^*(Q\times \mathbb{R})$, but, instead of referring to cosymplectic and symplectic structures on these spaces, we consider $T^*Q\times \mathbb{R}$ as a LCC manifold equipped with the differential forms $(\eta_{\Theta},\Omega_{\Theta})$ in \eqref{LCC-can} and $T^*(Q\times \mathbb{R})$ as a LCS manifold equipped with $\overline{\Omega}$ in   \eqref{LCS-Cos}.

To compute a geometric HJ theory for this realization we first plot the following commutative diagram referring to the fibrations in \eqref{pi-tau} and recall the sections and projections, as well as vector fields. 
 \begin{equation}\label{geomdiagram}
\xymatrix{ T^*(Q\times\mathbb{R})\ar[dr]^{\pi}
\ar[ddd]_{\pi_{Q\times \mathbb{R}}} \ar[rrr]^{X_{H^s}}&   & &TT^{*}(Q\times\mathbb{R})\ar[ddd]^{T\pi_{Q\times\mathbb{R}}}|<<<<<<<<<<<<{\hole} \ar[dr]^{T\pi}\\
  & T^* Q\times\mathbb{R}\ar[ddl]^{\tau}\ar[rrr]^{E_H\qquad \qquad}   & & &T(T^* Q\times\mathbb{R})\\\\
 Q\times\mathbb{R}\ar@/^3pc/[uuu]^{\tilde{\gamma}}\ar@/^1pc/[uur]^\gamma\ar[rrr]^{X_{H^s}^{\tilde{\gamma}}}_{E_H^{\gamma}}&  & & T(Q\times\mathbb{R})\ar[uur]_{T(\pi\circ \tilde{\gamma})}}
\end{equation}
In view of Diagram \ref{geomdiagram}, we have 
a section $\gamma$ of the fibration $\tau$ from $T^*Q\times \mathbb{R}$ to $Q\times \mathbb{R}$, and a section $\tilde{\gamma}$ of the cotangent bundle projection  $\pi_{Q\times \mathbb{R}}$ defined on $T^*(Q\times \mathbb{R})$. We recall the characterization of these sections as exhibited in \eqref{gamma-tilde}. Still we take the commutation of these sections under the projection  $\pi$ as depicted in \eqref{tildegamma-gamma}. 

For a time-dependent Hamiltonian function $H$ on the extended cotangent bundle $T^*Q\times \mathbb{R}$, the LCC evolution vector field $E_H$ is defined through \eqref{evo-LCC}. See that 
it is plotted in Diagram \ref{geomdiagram} as an arrow from $T^*Q\times \mathbb{R}$ to its tangent bundle $T(T^*Q\times \mathbb{R})$. The LCS Hamiltonian vector field for the extended Hamiltonian function $H^s$ is computed through \eqref{semiglobal}. Here, the extended Hamiltonian function $H^s$  is the one given in \eqref{H-s}. See that, in Diagram \ref{geomdiagram}, the Hamiltonian vector field $X_{H^s}$ is depicted as an arrow from $T^*(Q\times \mathbb{R})$ to its tangent bundle $TT^*(Q\times \mathbb{R})$. 
As in the case given in Section \ref{Sec-HJ-Cos}, we project the LCC evolution vector field $E_H$ to the manifold $Q\times \mathbb{R}$ by means of the projection $\gamma$, whereas we project the LCS Hamiltonian vector field $X_{H^s}$ to $Q\times \mathbb{R}$ by $\tilde{\gamma}$. So that we have
\begin{equation}\label{gamma-Ham---} 
X^{\Tilde{\gamma}}_{H^s} := T\pi_{Q\times\mathbb{R}} \circ X_{H^s} \circ \Tilde{\gamma},  \qquad 
E^{\gamma}_H  : =  T\tau \circ E_H \circ \gamma.  
\end{equation}
These projections are structurally the same as the ones in \eqref{gamma-Ham}. A direct computation reads that the projected vector fields are the same. 
They are pictured with a single arrow from $Q\times\mathbb{R}$ to its tangent bundle.

\begin{theorem}\label{gamma-rel-thm-LCCos}
Suppose that $\tilde{\gamma}$ is a LdR closed one-form on $Q\times \mathbb{R}$, that is $d_\Theta\tilde{\gamma}=0$. Following the notation presented in this section, the following conditions are equivalent:
\begin{enumerate}
\item $X^{\tilde{\gamma}}_{H^s}$ and $E_H$ are ${\gamma}$-related, that is,
\begin{equation}\label{LCCos-gamma-rel}
T\pi \circ T\tilde{\gamma} \circ X^{\tilde{\gamma}}_{H^s} = E_H \circ \pi \circ \tilde{\gamma}.
\end{equation}
\item The following identity holds:
\begin{equation}\label{HJ-eq_LCCos}
	d_\Psi(H^s \circ \tilde{\gamma}) \in \langle d_\Psi t \rangle,
\end{equation}
where $\langle d_\Psi t \rangle$ denotes the span space of $d_\Psi t$.
\end{enumerate}
\end{theorem}

\textbf{Proof.} Suppose that $X^{\tilde{\gamma}}_{H^s}$ and $E_H$ are $\gamma$-related. So, the pull-back $\gamma^*$ of $\gamma$ is distributive on the interior derivative, i.e.,
\begin{equation}\label{pi-gamma-related-iota-id2}
\tilde{\gamma}^*\pi^*\iota_{E_H}\Omega= \iota_{X^{\tilde{\gamma}}_{H^s}}\tilde{\gamma}^*\pi^*\Omega.
\end{equation}
The left hand side of \eqref{HJ-eq_LCCos} leads us to the following computation.
\begin{equation*}
\begin{split}
d_\Psi(H^s\circ \tilde{\gamma}) &= d_\Psi\tilde{\gamma}^*H^s = \tilde{\gamma}^* d_\Theta H^s \\
&= \tilde{\gamma}^*d_\Theta (\pi^*H + s) \\
&=  \tilde{\gamma}^*\pi^*d_\Theta H + d_\Psi\tilde{\gamma}^* s   \\
&= \tilde{\gamma}^*\pi^*\big(\iota_{E_H}\Omega +\langle d_\Theta H,\mathcal{R} \rangle\eta \big) + d_\Psi \bar{\bar{\gamma}} \\
&= \iota_{X^{\tilde{\gamma}}_{H^s}}\big(\tilde{\gamma}^*\pi^*\Omega 
\big)
+ \iota_{X^{\tilde{\gamma}}_{H^s}}  \tilde{\gamma}^*(\pi^*\eta \wedge u) - \iota_{X^{\tilde{\gamma}}_{H^s}} \tilde{\gamma}^*(\pi^*\eta \wedge u)  + \gamma^*(\langle d_\Theta H,\mathcal{R} \rangle\eta ) + d_\Psi \bar{\bar{\gamma}} \\
&= \iota_{X^{\tilde{\gamma}}_{H^s}}\tilde{\gamma}^*\overline{\Omega} - (\iota_{X^{\tilde{\gamma}}_{H^s}}\tilde{\gamma}^*\pi^*\eta) \tilde{\gamma}^*u + (\iota_{X^{\tilde{\gamma}}_{H^s}}\tilde{\gamma}^*u)\tilde{\gamma}^*\pi^*\eta + \gamma^*(\langle d_\Theta H,\mathcal{R} \rangle\eta ) + d_\Psi \bar{\bar{\gamma}} \\
&= 0 -(\tilde{\gamma}^*\pi^*\iota_{E_H}\eta) d_\Psi\bar{\bar{\gamma}} + (\iota_{X^{\tilde{\gamma}}_{H^s}}d_\Psi \bar{\bar{\gamma}})\gamma^*\eta + \gamma^*(\langle d_\Theta H,\mathcal{R} \rangle\eta) + d_\Psi \bar{\bar{\gamma}} \\
&= -d_\Psi\bar{\bar{\gamma}} +(\iota_{X^{\tilde{\gamma}}_{H^s}}d_\Psi \bar{\bar{\gamma}})\gamma^*\eta + \gamma^*(\langle d_\Theta H,\mathcal{R} \rangle\eta)+d_\Psi\bar{\bar{\gamma}} \\
&= \big(\iota_{X^{\tilde{\gamma}}_{H^s}}d_\Psi \bar{\bar{\gamma}} + \gamma^* \langle d_\Theta H,\mathcal{R} \rangle\big)\gamma^*\eta.
\end{split}
\end{equation*}
In the first line we have used the fact that the pull-back operation and the Lichnerowicz-deRham differential commute. In the second line, we have employed the definition of the Hamiltonian function $H^s$ in \eqref{H-s}. In the fourth line, we have used the definition \eqref{evo-LCC} of the evolution vector field $E_H$ and the explicit formulation \eqref{gamma-tilde} of the section $\tilde{\gamma}$. In the fifth  line, we have applied the identity \eqref{pi-gamma-related-iota-id2}, then we have added and subtracted the term $\iota_{X^{\tilde{\gamma}}_{H^s}}\tilde{\gamma}^*\pi^*\eta\wedge u$ in order to obtain the LCS two-form $\overline{\Omega}$ given in \eqref{Overline-Omega}. In the sixth line, we have used the distribution property of the pull-back operation $\tilde{\gamma}^*$ on wedge products, then we have applied the generalized Leibniz property of the interior derivative. In the seventh line, referring to the LdR closedness of $\tilde{\gamma}$, we have argued that $\tilde{\gamma}^*\overline{\Omega}$ vanishes identically. Further, we have applied the identity \eqref{pi-gamma-related-iota-id2} and have used the definition of $\tilde{\gamma}$. By the definition of $E_H$, the first and the last terms in the eight line cancel each other. The distributive character of the pull-back operation $\gamma^*$ provides the result in the last line. Since $\gamma^*\eta$ is indeed $d_\theta t$, we have arrived at the desired result. The converse statement can be proved through direct computation in local coordinates. $\blacksquare$

Notice that this theorem reduces to the Hamilton-Jacobi Theorem in cosymplectic manifolds \ref{gamma-rel-thm-Cos}. In this respect, we can argue that Theorem \ref{gamma-rel-thm-LCCos} is a generalization of Theorem \ref{gamma-rel-thm-Cos}.

\textbf{Hamilton-Poincar\'{e} Realization.} One can construct an alternative LCC structure on the extended cotangent bundle through the Hamilton-Poincar\'e one and two forms
\begin{equation}
\Theta_H = \Theta_Q - H d_\Theta t, \qquad \Omega_H = \Omega_{2\Theta} + d_\Theta H \wedge d_\Theta t,
\end{equation}
respectively. Here, the Hamiltonian-Poincar\'e two form can also be given by $\Omega_H = -d_{2\Theta}\Theta_H$. One can simply see that $(T^*Q\times\mathbb{R}, \eta,\Omega_H,\Theta)$ is a LCC manifold. We define the Reeb field $\mathcal{R}_H$ for this manifold through the following identities
\begin{equation}\label{reeblcc-HP}
\iota_{\mathcal{R}_H}\eta=1,\qquad \iota_{\mathcal{R}_H}\Omega_H=0.
\end{equation}
We can show that $R_H$ and the evolution vector field $E_H$ are the same by comparing the equations \eqref{reeblcc-HP} and \eqref{evo-LCC}. Clearly, the first identities in the definitions  are the same. Now let us evaluate the left hand side of the second equation in \eqref{reeblcc-HP}.
\begin{equation}
\begin{split}
\iota_{\mathcal{R}_H}\Omega_H &= \iota_{\mathcal{R}_H}\Omega_{2\Theta} + \iota_{\mathcal{R}_H}(d_\Theta H \wedge d_\Theta t) = \iota_{\mathcal{R}_H}\Omega_{2\Theta} + (\iota_{\mathcal{R}_H}d_\Theta H)d_\Theta t - (\iota_{\mathcal{R}_H}d_\Theta t)d_\Theta H \\
& = \iota_{\mathcal{R}_H}\Omega_{2\Theta} - d_\Theta H + \langle d_\Theta H, \mathcal{R}_H \rangle \eta.
\end{split}
\end{equation}
So, the second identity in \eqref{reeblcc-HP} yields
\begin{equation} \label{R-H-LCC}
\iota_{\mathcal{R}_H}\Omega_{2\Theta} = d_\Theta H - \langle d_\Theta H, \mathcal{R}_H \rangle \eta.
\end{equation}
Furthermore, the fact that $E_H = \mathcal{R} + X_H$ and the energy conservation allow us to write the second identity in \eqref{evo-LCC} as
\begin{equation} \label{E-H-LCC}
\iota_{E_H}\Omega_{2\Theta} = d_\Theta H - \langle d_\Theta H, E_H \rangle \eta.
\end{equation}
So the equations coincide, that is, we get $E_H=\mathcal{R}_H$.

An alternative construction of the Hamilton-Poincar\'{e} formulation can be performed by referring to the local cosymplectic formalism. To have this, suppose that $V_\alpha\subset T^*Q$ is an arbitrary open neighborhood admitting the Liouville one-form $\theta^\alpha_Q$ and the symplectic two-form $\omega^\alpha_Q = -d\theta^\alpha_Q$. We can construct an alternative cosymplectic structure on $U_\alpha = V_\alpha \times \mathbb{R}$ through the Hamilton-Poincar\'e one and two forms
\begin{equation}
\Theta_{H_\alpha} = \Theta^\alpha_Q-H_\alpha dt_\alpha,\qquad \Omega_{H_\alpha}=\Omega^\alpha _Q+dH_\alpha \wedge dt_\alpha,
\end{equation}
respectively, where $\Omega_{H_\alpha}$ is defined as $\Omega_{H_\alpha}=-d\theta_{H_\alpha}$. To obtain a global formulation for the extended cotangent bundle $T^*Q \times \mathbb{R}$, we compute
\begin{equation*}
\begin{split}
\Omega_{H_\alpha} &= \Omega^\alpha _Q+dH_\alpha \wedge dt_\alpha = e^{-2\sigma_\alpha}\Omega_Q\vert_\alpha + d(e^{-\sigma_\alpha}H\vert_\alpha) \wedge d(e^{-\sigma_\alpha}t\vert_\alpha) \\
&= e^{-2\sigma_\alpha}\Omega_Q\vert_\alpha + \big( e^{-\sigma_\alpha}dH\vert_\alpha -e^{-\sigma_\alpha}H\vert_\alpha d\sigma_\alpha \big) \wedge \big( e^{-\sigma_\alpha}dt\vert_\alpha -e^{-\sigma_\alpha}t\vert_\alpha d\sigma_\alpha \big) \\
&= e^{-2\sigma_\alpha} \big( \Omega_Q\vert_\alpha + (dH\vert_\alpha -H\vert_\alpha d\sigma_\alpha) \wedge (dt\vert_\alpha -t\vert_\alpha d\sigma_\alpha) \big)
\end{split}
\end{equation*}
The local two-forms
\begin{equation}
\Omega_H\vert_\alpha = e^{2\sigma_\alpha}\Omega_{H_\alpha}
\end{equation}
can be glued up to a global two-form $\Omega_H$.

\textbf{LdR Differential of a Function on $Q\times \mathbb{R}$.}
For a smooth function $F$ on $Q\times \mathbb{R}$, the LdR differential is one-form on $Q\times \mathbb{R}$. So that we can write it in the basis of $(d_\Psi q^i,d_\Psi t)$  as follows
\begin{equation}\label{d-phi-F}
d_\Psi F= F_{:i}  d_\Psi q^i + F_{:t} d_\Psi t.
\end{equation}
Here, $F_{:i} $ is the coefficient function of the basis element $d_\Psi q^i$ whereas $ F_{:t}$ is the coefficient function for the basis $d_\Psi t$. 
We determine the coefficients functions as 
\begin{equation}\label{coeff-F}
\begin{split}
F_{:i} & = F_{,i} + \frac{\psi_i}{1-q^k\psi_k - t\zeta}   \big(q^jF_{,j} + t F_{,t} -F\big)   , \\
F_{:t} & = F_{,t} + \frac{\zeta}{1-q^k\psi_k - t\zeta} \big(q^jF_{,j} + tF_{,t} -F \big) .
\end{split}
\end{equation}
Let us emphasize two subcases for these local formulas. One is the case $\Psi$ has no components along $dt$ that is $\zeta$ is zero. In this case, the coefficient functions are computed to be
\begin{equation}\label{coeff-F-zeta}
\begin{split}
F_{:i}  = F_{,i} + \frac{\psi_i}{1-q^k\psi_k }   \big(q^jF_{,j} + t F_{,t} -F\big), \qquad F_{:t}   = F_{,t}  .
\end{split}
\end{equation}
On the other extreme case, we choose that $\Psi$ has no components along $dq^i$, that is, $\psi_i$ is zero for all $i$. Then we have 
\begin{equation}\label{coeff-F-psi} 
F_{:i}  = F_{,i}, \qquad 
F_{:t}  = F_{,t} + \frac{\zeta}{1 - t\zeta} \big(q^jF_{,j}+ t F_{,t} -F \big) .
\end{equation}
Notice that, if $\Psi$ is identically zero then the LdR differential reduced to the deRham exterior derivative. In this case, the coefficient functions turn out to be the partial derivatives of the function $F$.  That is if both $\zeta$ and $\psi_i$ are zero, then $F_{:i}$ turns out to be the partial derivative of $F$ with respect to $q^i$ whereas $F_{:t} $ reduces to the partial derivative of $F$ with respect to $t$ as expected. See also that, if $F$ is a homogeneous function of order $1$ then one has that $q^j F_{,j} + t F_{,t} =F$. In this case, $F_{:i}  = F_{,i}$ and $F_{:t}  = F_{,t}$ for any closed one-form $\Psi$. 

\textbf{LdR Differential of a One-form on $Q\times \mathbb{R}$.} We have examined the local formulation of LdR differential in Appendix \ref{Sec-LdR}. Let us now carry this discussion to the extended configuration space as follows.  
Assume a local coordinate chart $(q^i)$ on a manifold $Q$ and fix a closed one-form $\Psi$ in the form \eqref{Lee-Psi} on the extended manifold $Q\times \mathbb{R}$.  
Let us observe that the set $(d_\Psi q^i,d_\Psi t)$ determines a basis for the module of one-form sections on $Q\times \mathbb{R}$. So that, we can write a differential one-form $\upsilon$ as 
\begin{equation}\label{upsilon}
\begin{split}
  \upsilon(q,t)&=\hat{\upsilon} _i (q,t) ~ d_\Psi q^i+ \hat{\hat{\upsilon}} (q,t)~  d_\Psi t, \\
  \upsilon(q,t)&=\bar{\upsilon} _i (q,t) ~ d  q^i+ \bar{\bar{\upsilon}} (q,t)~  d  t, 
\end{split}
\end{equation}
A direct calculation proves the following relations between the coefficient functions $(\hat{\upsilon} _i,\hat{\hat{\upsilon}})$ with respect to the basis $(d_\Psi q^i,d_\Psi t)$ and the coefficient functions $(\bar{\upsilon} _i,\bar{\bar{\upsilon}}$ with respect to the basis $(dq^i,dt)$ as follows 
\begin{equation}
\begin{split}
     \bar{\upsilon}_{i} & =\hat{\upsilon}_{i}-\psi _{i}(q^{l}\hat{\upsilon}_{l}+t\hat{\hat{\upsilon}}), \qquad \bar{\bar{\upsilon}} =\hat{\hat{\upsilon}} - \zeta (q^{i}\hat{\upsilon}_{i} 
+t  \hat{\hat{\upsilon}}), 
\\
\hat{\upsilon}_{i} &=\bar{\upsilon}_{i}+\frac{q^{j}\bar{\upsilon}_{j}+\bar{\bar{\upsilon}} t}{%
1-q^{i}\psi _{i}-t\zeta }\psi _{i} ,\qquad 
\hat{\hat{\upsilon}} =\bar{\bar{\upsilon}} +\frac{q^{j}\bar{\upsilon}_{j}+\bar{\bar{\upsilon}} t}{1-q^{i}\psi
_{i}-t\zeta }\zeta    
\end{split}
\end{equation}
provided that $1-q^{i}\psi
_{i}-t\zeta $ is not zero. Let us compute the conditions of being LdR closed for both of the local forms of the one-form $\upsilon$ given in \eqref{upsilon}. 
For the former local realization of $\upsilon$ determined through the the LdR basis $\{d_\Psi q^i,d_\Psi t\}$, we compute that $d_\Psi \upsilon$ is zero if and only if 
\begin{equation}
\begin{split}
\hat{\upsilon}_{i,j}+\hat{\upsilon}_{k,i}q^{k}\psi _{j}+t\hat{\hat{\upsilon}}
_{,i}\psi _{j}  &= \hat{\upsilon}_{j,i}+\hat{\upsilon}_{k,j}q^{k}\psi _{i}+t%
\hat{\hat{\upsilon}}_{,j}\psi _{i}
\\
\hat{\hat{\upsilon}}_{,i} +\hat{\hat{\upsilon}}_{,t}t\psi
_{i}-t\zeta \hat{\hat{\upsilon}}_{,i} &=\hat{\upsilon}_{i,t}+q^{j}\zeta \hat{\upsilon}_{j,i}-\hat{\upsilon}_{k,t}q^{k}\psi _{i}.
\end{split}
\end{equation}
For the latter realization in \eqref{upsilon} of $\upsilon$ done for the case of the basis $\{dq^i,dt\}$, we compute the conditions of being LdR closed as 
\begin{equation}\label{LdR-closed-dup-}
\bar{\upsilon}_{i,j}+\psi _{i}\bar{\upsilon}_{j} =\bar{\upsilon}%
_{j,i}+\psi _{j}\bar{\upsilon}_{i},\qquad 
\bar{\upsilon}_{i,t}+\bar{\bar{\upsilon}}\psi _{i}=\bar{\bar{\upsilon}}_{,i}+\zeta 
\bar{\upsilon}_{i}.
\end{equation}

\textbf{HJ Theorem for Evolution Vector Field.}
To have this, we consider a smooth function $S=S(q,t)$ defined on the extended configuration space $Q\times \mathbb{R}$. Then its exterior derivative determines the following section 
\begin{equation}\label{tilde-dS}
\tilde{\gamma}(q,t)=dS(q,t)=(q^i,S_{,i}(q,t),t,S_{,t}(q,t))\in T^*Q\times T^*\mathbb{R}
\end{equation}
given in terms of the canonical coordinates $(q^i,p_i,t,s)$. We wish that $dS$ is LdR closed. In the light of the local restrictions exhibited in \eqref{LdR-closed-dup-}, $d_\Psi dS$ is zero if and only if
\begin{equation}\label{dS-LdR-closed}
\psi _{i}S_{,j}  = \psi _{j}S_{,i}, \qquad S_{,t}\psi _{i} =\zeta S_{,i}.
\end{equation}
Recall  the second item in Theorem \ref{gamma-rel-thm-LCC} and then examine it by taking the section $\tilde{\gamma}$ as in the form of \eqref{tilde-dS}. Then we compute
\begin{equation}
\begin{split}
d_\Psi (H^s\circ \tilde{\gamma})&=d_\Psi (\pi^*(H(q^j,S_{,j}(q,t),t) + S_{,t})
\\
&= d_\Psi \pi^*(H(q^j,S_{,j}(q,t),t))+  d_\Psi S_{,t}
\\
&= \pi^*d_\Theta(H(q^j,S_{,j}(q,t),t))+ \big((S_{,t})_{:i}d_\Psi q^i + (S_{,t})_{:t}d_\Psi t \big)
\\
&=\big(H(q^j,S_{,j}(q,t),t)_{:i}+(S_{,t})_{:i}\big) d_\Psi q^i+  
\big(
H(q^j,S_{,j}(q,t),t)_{:t}+(S_{,t})_{:t}
\big)d_\Psi t.
\end{split}
\end{equation}
If we insist that this should be parallel to $d_\Psi t$, then we can claim that
\begin{equation}\label{HJ-111}
H(q^j,S_{,j}(q,t),t)_{:i}+(S_{,t})_{:i} =0.
\end{equation}
This is the Hamilton-Jacobi equation for the LCC Hamiltonian dynamics. We are ready now to state the following version of the geometric Hamilton-Jacobi theorem for LCC framework. 

\begin{theorem}\label{gamma-rel-thm-LCC} Consider a LCC manifold $(T^*Q\times \mathbb{R},\eta,\Omega,\Theta)$, and a section 
\begin{equation}
\gamma(q,t)=(q^i,S_{,i}(q,t),t)
\end{equation}
of the fibration $\tau$. Assume that $dS$ is LdR closed that is the local conditions in \eqref{dS-LdR-closed} are satisfied. 
Then the following conditions are equivalent:
\begin{enumerate}
\item The vector fields $E_H$ and $E_H^{\gamma}$ are $\gamma$-related.
\item The following identity holds:
\begin{equation}\label{HJ-LCC-eq}
H(q^j,S_{,j}(q,t),t)_{:i}+(S_{,t})_{:i} =0.
\end{equation}
\end{enumerate}
\end{theorem}
Notice that, the second condition determines a Hamilton-Jacobi equation for a locally conformal cosymplectic Hamiltonian flow. To have it more explicitly, one can refer to \eqref{coeff-F}. For the case of LCC structures those generated by Lee fields independent of time variable (that is $\zeta=0$ case), then the Hamilton-Jacobi equation reduces to a more simple form and can be written directly to \eqref{coeff-F-zeta}. On the other hand, for the case of LCC structures those generated by Lee fields independent of space variable (that is $\psi_i=0$ case), then the Hamilton-Jacobi equation reduces to a more simple form and can be written directly to \eqref{coeff-F-psi}. If both of these coefficient functions are zero then the LCC geometry becomes cosymplectic geometry and, accordingly, Theorem \eqref{gamma-rel-thm-LCC} particularly realizes Theorem \ref{thm-Cristina-cos} given for pure cosymplectic formalism. As in the case of pure cosymplectic case, in Theorem \ref{gamma-rel-thm-LCC}, we may substitute the evolution vector field $E_H$ and its projection $E^{\gamma}_H$ with the Reeb field $\mathcal{R}_H$ and $\mathcal{R}^{\gamma}_H$, respectively. As discussed previous, this can be directly justified by comparing \eqref{R-H-LCC} and \eqref{E-H-LCC}.

We would like to remark the extreme case $\Psi = \zeta(t)$ forces $S_{,q}$ to be zero, that is, $S=S(t)$. As a result, we cannot obtain a PDE from the identity \eqref{HJ-LCC-eq}.

\subsection{An Illustration: Damped Harmonic Oscillator}

Consider the $3$-dimensional space with  coordinates $(q,p,t)$. Here we are interested in a damped harmonic oscillator that has a time-dependent mass of the form $m = m_0e
^{\Gamma t}$ accreting with time in order to mimic exponential energy dissipation, see \cite{GaseGracMuno20}. For its cosymplectic dynamic one can check some results in \cite{esen2022reviewing}. The dynamic is generated by the Hamiltonian function 
\begin{equation}\label{Ham-ex}
H(q,p,t)=\frac{p^2}{2m}e^{-\Gamma t}+\frac{m}{2}e^{\Gamma t}q^2.
\end{equation}
We assume that the dynamic here is local and the LCC framework is determined through the Lee-form \eqref{Lee-LCCos} which is independent of the position coordinate
\begin{equation}\label{tht}
\Theta(q)=\Psi(q)=\psi(q)dq.
\end{equation}
The form \eqref{tht} is evidently a closed form. According to \eqref{LCC-can}, we notice the following one-form and two-form
   \begin{equation} \label{forms-ex}
    \Omega_{2\Theta} = dq \wedge dp,  \qquad    \eta = d_{\Theta}t = dt - t\psi dq, 
     \end{equation}
     respectively. Through \eqref{local-Reeb-LCC} we compute the Reeb field on LCC as: 
\begin{equation}
\mathcal{R} = \frac{\partial}{\partial t}.
\end{equation}  
Recall that the closed one-form $\Psi$ determines a new basis for the space of one-form section 
   \begin{equation}  
d_\Psi q = (1-q\psi)dq, \qquad d_\Psi t = dt - t\psi dq. 
     \end{equation}
With respect to this basis and referring to the local expression in \eqref{d-phi-F}, the LdR differential of a function $F = F(q,t)$ admits the following coefficient functions 
\begin{equation}\label{coeff-F-psi-ex}
F_{:q}  = F_{,q} + \frac{\psi}{1-q\psi }   \big(qF_{,q} + t F_{,t} -F\big), \qquad F_{:t}   = F_{,t} .
\end{equation}

Let us examine the gluing problem for the dynamics to this physical system. To have this we recall the evolution vector field for LCC geometry given in \eqref{E-H-LCC-local} then,  compute the evolution vector field as 
\begin{equation}
\begin{split}
    E_H = \frac{p}{m}e^{-\Gamma t} \frac{\partial}{\partial q} &+ \left(-me^{\Gamma t}q + \psi(\frac{p^2}{2m}e^{-\Gamma t}+\frac{m}{2}e^{\Gamma t}q^2) -t\psi\Gamma\left(\frac{m}{2}e^{\Gamma t}q^2 - \frac{p^2}{2m}e^{-\Gamma t}\right)\right) \frac{\partial}{\partial p} \\
    &+ \left( 1 + t\psi \frac{p}{m} e^{\Gamma t} \right) \frac{\partial}{\partial t}.
\end{split}
\end{equation}
Therefore, the Hamilton equations are
\begin{align}
    \dot{q} &= \frac{p}{m}e^{-\Gamma t}, \nonumber \\
    \dot{p} &= -me^{\Gamma t}q + \psi\left(\frac{p^2}{2m}e^{-\Gamma t}+\frac{m}{2}e^{\Gamma t}q^2 \right) -t\psi\Gamma\left(\frac{m}{2}e^{\Gamma t}q^2 - \frac{p^2}{2m}e^{-\Gamma t}\right), \\ 
    \dot{t} & = 1 + t\psi \frac{p}{m} e^{\Gamma t}.\nonumber
 \end{align}   
Now, let us choose a section as
\begin{equation}
   \gamma(q,t)=(q,S_{,q}(q,t),t),
\end{equation}
and consider a LdR closed one-form
\begin{equation} 
\tilde{\gamma}(q,t)=dS(q,t)=(q,S_{,q}(q,t),t,S_{,t}(q,t)),\qquad \tilde{\gamma}(q,t)= S_{,q} dq + S_{,t} dt,
\end{equation}
where $S=S(q,t)$ is a smooth function. The LdR-closedness of $\tilde{\gamma}(q,t)$ imposes $S_{,t} = 0$, that is, $S=S(q)$.
Then the projected vector field $E_H^\gamma$ determined in \eqref{gamma-Ham---} is computed to be 
\begin{equation} 
E_H^\gamma  = \frac{S_{,q}}{m}e^{-\Gamma t} \frac{\partial}{\partial q} + \left( 1 + t\psi \frac{p}{m} e^{\Gamma t} \right) \frac{\partial}{\partial t}.
\end{equation}
Furthermore, the HJ equation \eqref{HJ-LCC-eq} for the LCC evolution dynamics is
\begin{equation}
\Big(\frac{S_{,q}^2}{2m}e^{-\Gamma t}+\frac{m}{2}e^{\Gamma t}q^2 \Big) _{:q}+ (S_{,t})_{:q}=0.
\end{equation}
We use the fact that $S_{,t} = 0$ and the first equation in \eqref{coeff-F-psi-ex} to write the HJ equation as
\begin{equation}
    \frac{1 }{1-q\psi}\left(\frac{1}{m}S_{,q}(S_{,q})_{,q}e^{-\Gamma t} + me^{\Gamma t}q\right) = (1 + \Gamma t)\frac{S_{,q}^2}{2m}e^{-\Gamma t} + (1 - \Gamma t)\frac{m}{2}e^{\Gamma t}q^2
\end{equation}

\section{Conclusion and Future Works}

In this work, we have addressed the gluing problem of local time-dependent Hamiltonian dynamics. This is achieved in the realm of locally conformal cosymplectic manifolds. The symplectization and Darboux coordinates of LCC manifolds are obtained (c.f. Section \ref{Sec-LCS-LCCos}). The Jacobi structure (c.f. Section \ref{sec-Jac-LCC}) and the Lie algebroid structure (c.f. Section \ref{Sec-LCC-Algebroid}) of locally conformally cosymplectic manifolds are investigated.  Dynamically, we have obtained a global picture of local time-dependent Hamiltonian systems (c.f. Section \ref{Sec-Dyn-LCCos}) and we have presented their associated geometric Hamilton-Jacobi theories (c.f. Section \ref{Sec-HJ-LCC}). There are some future directions we would like to continue working on. Let us list them and comment on them one by one. 

\begin{itemize}

\item \textit{A Generalization Through Lagrangian/Legendrian Submanifold Realizations.} One of the generalizations of the classical Hamilton-Jacobi equation (in the symplectic setting) is based on the replacement of a generating function $S$ by a Lagrangian submanifold, see \cite{BarbLeonDieg,BeneTulc80}. We wish to focus on this generalization in the future. In this case, we shall search for all possible Lagrangian/Legendrian submanifolds of LCS and LCC manifolds including the non-horizontal ones (those which cannot be written in terms of generating functions, not even locally), see \cite{Benenti-book,SniaTulc}. This may permit us to arrive at complete solutions of the Hamilton-Jacobi equations given in the present work. 

\item \textit{The Reduction of LCC Manifolds Under Symmetry.}
The reduction of Hamiltonian systems under Lie group symmetries, aka. Marsden-Weinstein reduction \cite{marsden74,Marsden1999} has been generalized in many different geometric scenarios. 
For example, the reduction of LCS manifolds has been studied in \cite{Stanciu19,Stan22}, the reduction for Jacobi manifolds has been addressed in \cite{IborLeonMarm97} and the reduction of cosymplectic manifolds has been obtained in \cite{Albe87}. 
We would like to address the reduction of LCC manifolds under a Lie group action in an upcoming work. 

\item \textit{Implicit Cosymplectic Hamiltonian Dynamics.} The dynamical systems derived in this paper are all described in terms of explicit differential equations. In a recent work, we have proposed the Hamilton-Jacobi Theorem for implicit Hamiltonian dynamics in the symplectic formalism \cite{EsenLeonSar1,EsenLeonSar2}. As future work, we plan on constructing an implicit HJ theory, as well as implicit Hamiltonian dynamics for cosymplectic geometry.

\item \textit{Locally Conformal $k$-cosymplectic Geometry.} Hamiltonian dynamics on $k$-cosymplectic manifolds was proposed in \cite{LeonMaeriOubiRodrSalg98}, and the geometric HJ theorem on such manifolds is available in  \cite{LeonVila14}, whilst the $k$-symplectic approach is in \cite{LeonMartMarrSalgVila10}. One can find further notes on these two theories in  \cite{LeonSalgVila16}. On the other hand, the geometric HJ for locally conformal $k$-symplectic has recently been established in \cite{EsLeSaZa-k-sympl}. Locally conformal $k$-cosymplectic formalism is still missing in the literature, so we hope to address it in our upcoming work.

\item \textit{Tulczyjew's Triplet for LCS and LCC Formalisms.} The Tulczyjew's triplet \cite{tul77} is a geometric construction that permits the existence of a  Legendre transformation between Hamiltonian and Lagrangian dynamics even for degenerate cases. This triplet has been generalized in different geometric frameworks, as it is the case of time-dependent dynamics in \cite{LeonMarr93}, $k$-cosymplectic framework in \cite{ReyRomNarSalVil12} and recently, on contact manifolds  \cite{esen2021contact}.  We wish to construct Tulczyjew triples for LCS and LCC manifolds.

\item \textit{Analysis on Lie Algebroid Realizations of LCC Manifolds.} The cotangent bundle of a Poisson manifold admits a Lie algebroid structure \cite{CostDazoWein87}. This has been generalized to Jacobi manifolds in \cite{LeonMarr97}, see also \cite{LeonMarrJuan97}. We would like to examine the Lie algebroid realizations obtained in this work following the approach presented in \cite{LeonMarr97,LeonMarrJuan97}.

\end{itemize}

\section*{Acknowledgements}

BA and OE gratefully acknowledge Prof. Serkan Sütlü for discussions on Lie algebroids. Manuel de León acknowledges financial support from the Spanish Ministry of Science and Innovation (MICINN), under grants PID2019-106715GB-C21 and “Severo Ochoa Programme for Centres of Excellence in R\&D” (CEX2019-000904-S). CS acknowledges project “Teoría de aproximación constructiva y aplicaciones” (TACA-ETSII).

\newpage
\appendix 
\section{Appendix}
\subsection{The Lichnerowicz-deRham Differential}\label{Sec-LdR} 

Consider now an arbitrary manifold $M$. Fix a closed one-form $\theta$.  The Lichnerowicz-deRham differential is defined as
\begin{equation} \label{LdR-Diff}
d_\theta: \Lambda^k(M) \rightarrow  \Lambda^{k+1}(M) : \beta \mapsto d\beta-\theta\wedge\beta,
\end{equation}
where $d$ denotes the exterior (deRham) derivative \cite{GuLi84}. See that the Lichnerowicz-deRham (abbreviated as LdR) differential is an idempotent operator, i.e.,  $d_{\theta}^2=0$. See, for example, \cite{ChanMurp19} for more properties and some details on cohomological discussions.  

\textbf{LdR Differential of a Function.} We assume a local chart $\{x^a\}$ on $M$ and consider the closed one-form $\theta=\theta_a dx^a$. 
See that, the set $\{d_\theta x^a\}$ determines a basis for the space of sections if $1-x^a\theta_a$ does not vanish for any $x$ in $M$. 
In this case, for a smooth function $F=F(x^a)$  on $M$, we can write a LdR exact one-form $d_\theta F$ in terms of the basis $\{d_\theta x^a\}$ as 
\begin{equation}\label{theta-decom}
d_\theta F = F_{:a} d_\theta x^a.
\end{equation}
Here, $F_{:a}$ are the coefficient functions computed to be
\begin{equation} 
F_{:a}   = F_{,a} + \frac{\theta_a }{1-x^c\theta_c  } \left( x^bF_{,b}  -F \right).
\end{equation}
 Notice that, if $\theta$ identically vanishes, the LdR reduces to the deRham exterior derivative and $F_{:a}$ turns out to be the partial derivative $F_{,a}$ of $F$ with respect to the coordinate $x^a$. See also that if $F$ is a homogeneous function of order $1$, thanks to the Euler's theorem $x^b F_{,b}= F$, then $F_{:a}   = F_{,a} $.
 
 \textbf{LdR Differential of a One-form.}
As indicated in the previous paragraph, for a coordinate chart $\{x^a\}$ on $M$, we can introduce two different bases for the module of one-form sections, namely $\{dx^a\}$ and $\{d_\theta x^a\}$. Accordingly, we can write the local realization of one-form $\upsilon$ as
\begin{equation} \label{upps}
\upsilon =\bar{\upsilon}_{c}dx^{c}, \qquad \upsilon =\hat{\upsilon}_{c}d_{\theta }x^{c}, \qquad \bar{\upsilon}_{a}= \hat{\upsilon}_{a}-\hat{\upsilon}_{c}x^{c}\theta _{a}.
\end{equation} 
Following the order of the local realization in \eqref{upps}, the LdR-closedness of  $\upsilon$ dictates the following symmetries for the one-form $\upsilon$ 
\begin{equation}\label{dLR-closed}
    \bar{\upsilon}_{c,a}+\theta _{c}\bar{\upsilon}_{a} =\bar{\upsilon}%
_{c,a}+\theta _{a}\bar{\upsilon}_{c},\qquad 
\hat{\upsilon}_{c,b}+\hat{\upsilon}_{a,c}\theta _{b}x^{a} =\hat{\upsilon}%
_{b,c}+\hat{\upsilon}_{a,b}\theta _{c}x^{a},
\end{equation}
respectively. For some smooth function $F$ on $M$, we can state that natural solutions to these conditions \eqref{dLR-closed} are $\bar{\upsilon}_c=F_{:,c}$ and $\hat{\upsilon}_c=F_{:c}$, respectively.  Such a local solution always exists for LdR closed forms, and this result is called the Poincar\'{e} lemma for LdR differential. See that if $\theta=0$, then one arrives at the well-known Poincar\'{e} lemma. 
\subsection{Jacobi Manifolds}\label{Sec-Jac-M}

 A manifold $M$ equipped with a vector field $\mathcal{E}$ and a bivector field $\Lambda$ is a Jacobi manifold if 
\begin{equation}\label{ident-Jac}
        [\Lambda,\Lambda] = 2 \mathcal{E} \wedge \Lambda, \qquad 
       [\mathcal{E},\Lambda] = 0,
\end{equation}
    where $[\bullet,\bullet ]$ is the Schouten-Nijenhuis bracket, see, for example, \cite{Kirillov,Lichnerowicz-Poi,Lichnerowicz-Jacobi,Marle-Jacobi,VitaWade20}. We denote a Jacobi manifold by a triple $(M,\Lambda,\mathcal{E})$. Starting with a Jacobi manifold, one may define an antisymmetric bilinear bracket
   \begin{equation}\label{bra-Jac}
        \{H,F\} = \Lambda(dH, dF) + H\mathcal{E}(F) - F\mathcal{E} (H).
 \end{equation}
satisfying the Jacobi identity. Furthermore, it fulfills the so-called the weak Leibniz rule
\begin{equation}
        \operatorname{supp}(\{F,G\}) \subseteq \operatorname{supp} (F) \cap \text{supp} (G).
\end{equation}
This observation reads that the algebra is a local Lie algebra in the sense of Kirillov \cite{Kirillov}. The inverse of this assertion is also true, that is, a local Lie algebra determines a Jacobi structure. 

Now assume a differentiable map $\chi$ from a Jacobi manifold $\ (M_{1}, \Lambda_{1}, \mathcal{E}_{1} )$ to another Jacobi manifold $ (M_{2}, \Lambda_{2}, \mathcal{E}_{2} )$. This map is said to be a Jacobi map if it respects the induced Jacobi bracket, 
\begin{equation}
\{F \circ \chi, G \circ \chi\}_{1}=\{F, G\}_{2} \circ \chi
\end{equation}
for all smooth functions $F$ and $H$ on $M_{2}$. 

Consider a Jacobi manifold determined by the triplet $(M, \Lambda, \mathcal{E})$. 
We determine a homomorphism from the space of one-forms to the space of vector fields. It is determined by
\begin{equation}
\sharp_\Lambda:\Gamma^{1}(M) \longrightarrow \mathfrak{X}(M),\qquad \langle \sharp_\Lambda(\mu), \nu \rangle =\Lambda(\mu, \nu)
\end{equation}
for all $\mu$ and $\nu$ in $\Gamma^{1}(M)$. For a smooth  real-valued Hamiltonian function $H$, the  Hamiltonian vector field $X_{H}$ is defined by
\begin{equation}
X_{H}=\sharp_\Lambda(d H)+H \mathcal{E}.
\end{equation}
Note that, for the constant function $H=1$ the Hamiltonian vector field is $\mathcal{E}$. As proved in \cite{Lichnerowicz-Jacobi,Marle85}, the mapping taking a Hamiltonian function $H$ to the Hamiltonian vector field $X_H$ is a Lie algebra homomorphism satisfying
\begin{equation}
\left[X_{F}, X_{H}\right]=X_{\{F, H\}}.
\end{equation}

\subsection{Lie Algebroids}\label{Defn-Lie-Algebroid}

Given a manifold $M$, a Lie algebroid $\mathcal{A}$ over the base $M$ is a (real) vector bundle $\tau:\mathcal{A} \to M$, together with a map $a:\mathcal{A} \to TM$ of vector bundles, called the anchor map, and a Lie bracket $[\bullet,\bullet]$ (bilinear, antisymmetric, satisfying the Jacobi identity) on the space $\Upsilon(\mathcal{A})$ of sections, so that the induced $\mathcal{F}(M)$-module homomorphism   $a$ from $\Upsilon(\mathcal{A})$ to the space $\mathfrak{X}(M)$ of vector field on $M$ satisfies
\begin{equation}\label{oid-id}
[\xi, F\eta] = F[\xi,\eta]+\mathcal{L}_{a(\xi)}(F)\eta
\end{equation}
for any $\xi,\eta \in \Gamma(\mathcal{A})$, and any $F\in \mathcal{F}(M)$, where $\mathcal{L}_{a(\xi)}(F)$ stands for the Lie derivative of $ F$ in the direction of $a(\xi)$ in $TM$. See \cite{Mackenzie-book,Mackenzie-book-II,Para67}. Then, it follows that 
\begin{equation}
a([\xi,\eta]) = [a(\xi),a(\eta)]
\end{equation}
for any $\xi,\eta \in \Gamma(\mathcal{A})$. Accordingly, a Lie algebroid is denoted by a quintuple $(\mathcal{A},\tau,M,a,[\bullet,\bullet])$.

The cotangent bundle of a Poisson manifold admits a Lie algebroid structure \cite{BhasVisw88,CostDazoWein87}. For a Jacobi manifold $(M, \Lambda, \mathcal{E})$, the picture is as follows. Consider the extended cotangent bundle $T^{*} M \times \mathbb{R}$ as the total space of the first jet prolongation of the fibration $M\mapsto \mathbb{R}$, see \cite{Saunders-book}. The space of sections of $p_M:T^{*} M \times \mathbb{R}\mapsto M$ is the product space of one-forms and real valued functions $\Gamma^1(M)\times \mathcal{F}(M)$. So it consists of pairs $(\mu,F)$ where $\mu$ is a one-form and $F$ be a real valued function. A bracket on $\Gamma^1(M)\times \mathcal{F}(M)$ can be defined to be
\begin{equation}\label{Jac-Brack}
\begin{split} 
\{(\mu, F),(\nu, H)\}_J&=\Big(\mathcal{L}_{\sharp_\Lambda(\mu)} \nu-\mathcal{L}_{\sharp_\Lambda(\nu)} \mu-d(\Lambda(\mu, \nu))+F \mathcal{L}_{\mathcal{E}} \nu-H \mathcal{L}_{\mathcal{E}} \mu-\iota_{\mathcal{E}}(\mu \wedge \nu),\\
&\qquad  \mu(\sharp_\Lambda(\nu))+\sharp_\Lambda(\mu)(H)-\sharp_\Lambda(\nu)(F)+F \mathcal{E}(H)-H \mathcal{E}(F)\Big),
\end{split}
\end{equation}
where $\sharp_\Lambda$ is the musical mapping induced by the bivector field $\Lambda$. 
Referring to this bracket, the quintuple is $(T^{*} M \times \mathbb{R},p_M,M,\sharp_\Lambda+\mathcal{E},\{\bullet,\bullet\}_J)$, see \cite{KerbSoui93}. Here the anchor map is given by 
\begin{equation}
\sharp_\Lambda+\mathcal{E}:\Gamma^1(M)\times \mathcal{F}(M)\longrightarrow \mathfrak{X}(M),\qquad (\mu,F)\mapsto \sharp_\Lambda(\mu)+F\mathcal{E}. 
\end{equation}
We cite \cite{LeonMarr97,LeonMarrJuan97,LeonLopeMarrPadr03} for some further analysis on this Lie algebroid. If $\mathcal{E}$ is identically zero, then one can consider the cotangent bundle $T^*Q$ as the total space and the bracket \eqref{Jac-Brack} reduces to the algebra of one-forms on the Poisson manifold hence the Poisson Lie algebroid structure. We remark here that the Lie algebroid structures provided for LCS and LCC manifolds (those given in Section \ref{Sec-LCS-Algebroid} and Section \ref{Sec-LCC-Algebroid}, respectively) are not in the form presented here. We have employed the musical isomorphisms (considered to be the anchor maps) in both LCS and LCC geometries.

\bibliographystyle{abbrv}
\bibliography{references}

\end{document}